\documentclass[final,leqno,onefignum,onetabnum]{siamltex}
\pdfoutput=1        

\usepackage{algorithm}
\usepackage{algpseudocode}
\usepackage{pifont}
\usepackage{amsmath,subfigure,bm}
\usepackage{amssymb}
\usepackage{graphics}
\usepackage{graphicx}
\usepackage{color}          
\usepackage{cite}           

\usepackage{hyperref}       

\newcommand\Z{\mathbb{Z}}

\newcommand\R{\mathbb{R}}

\newcommand{\Tm}{\mathcal{T}}

\newcommand{\dt}{\Delta t}
\newcommand{\dx}{\Delta x}
\newcommand{\dy}{\Delta y}

\newcommand{\hgrid}{h}

\newcommand{\BigOh}{\mathcal{O}}

\newtheorem{rmk}{Remark}

\newcommand{\En}{\mathcal{E}}

\def\reals{\mathbb R}

\def\WS{{\mathcal W}}
\def\deg{{M_{\text{D}}}}
\def\degtwo{{M^2_{\text{D}}}}
\def\mleg{{M_{\text{L}}}}
\def\beq{\begin{equation}}
\def\eeq{\end{equation}}
\def\meq{{M_{\text{E}}}}
\def\bxi{\boldsymbol{\xi}}
\def\half{\frac{1}{2}}
\def\ubold{{\bf u}}

\title{A Simple and Effective High-Order Shock-Capturing Limiter for 
Discontinuous Galerkin Methods}

\author{Scott A. Moe\thanks{University of Washington, Department of Applied Mathematics,
Box 353925, Seattle, WA 98195, USA ({\tt smoe@uw.edu})}
\and James A. Rossmanith\thanks{Iowa State University, Department of Mathematics,
396 Carver Hall, Ames, IA 50011, USA ({\tt rossmani@iastate.edu})}
\and David C. Seal\thanks{Michigan State University, Department of Mathematics, 619 Red Cedar Road,
East Lansing, MI 48824, USA ({\tt seal@math.msu.edu})}}

\begin{document}
\maketitle
\slugger{sisc}{xxxx}{xx}{x}{x--x}

\begin{abstract}
The discontinuous Galerkin (DG) finite element method when applied to hyperbolic conservation laws requires the use of {\it shock-capturing limiters} in order to suppress unphysical oscillations near large solution gradients. In this work we develop a novel
shock-capturing limiter that combines key ideas from the limiter of Barth and Jespersen
[{\it AIAA-89-0366} (1989)] and the maximum principle 
preserving (MPP) framework of Zhang and Shu
[{\it Proc. R. Soc. A}, {467} (2011), pp. 2752--2776].
 The limiting strategy is based on traversing the mesh element-by-element in order to
 (1) find local upper and lower bounds on user-defined variables by
 sampling these variables on neighboring elements, and (2) to then enforce these local bounds by minimally damping the high-order corrections. The main advantages of this
 limiting strategy is that it is simple to implement, effective at shock capturing, and retains
high-order accuracy of the solution in smooth regimes.  
The resulting numerical scheme is applied to several standard numerical tests in both one 
and two-dimensions and on both Cartesian and unstructured grids. These
tests are used as benchmarks to verify and assess the accuracy and robustness
of the method.
\end{abstract}

\section{Introduction}
\subsection{Some background}
Godunov's theorem \cite{article:Go59} states that a {\it linear} numerical
method when applied to the one-dimensional, linear, advection equation:
\begin{gather}
\label{eqn:advection1d}
	\frac{\partial q}{\partial t} + \frac{\partial q}{\partial x} = 0, \quad q(t=0,x) = q_0(x),
\end{gather}
where $x \in \reals$, $t\in \reals_{\ge 0}$, and $q(t,x) : \reals_{\ge 0} \times \reals
\rightarrow \reals$, is {\it monotone} (i.e., does not generate new extrema) only if it is at most first-order accurate. In practice, schemes that are not monotone will exhibit large
spurious oscillations near discontinuities in the solution.
The analog of this result for more complicated hyperbolic equations (i.e., nonlinear systems in multiple dimensions) is that linear numerical methods with a formal accuracy higher than
first-order will create spurious oscillations that can lead to catastrophic numerical
instabilities (e.g., negative densities and/or pressures). 
As in the linear scalar case, these instabilities are especially pronounced if the
solution becomes discontinuous (i.e., shock waves and contact discontinuities).
The additional challenge with nonlinear hyperbolic equations is that 
generic initial conditions tend to produce shocks in finite time.

The principal remedy for these numerical instabilities in high-order methods 
is to introduce a {\it shock-capturing limiter}; such limiters render an otherwise 
linear method nonlinear, thereby circumventing the restrictions of 
Godunov's theorem. Typically, limiters are designed to execute two tasks:
(1) to monitor the solution quality (at least in some simplified sense); 
and, when necessary, (2) to reduce the influence of high-order corrections.
The development of limiters for hyperbolic equations was pioneered by 
Bram van Leer in a series of papers starting with \cite{article:vanleer73,van1974towards}. 
Additional fundamental contributions came from Harten \cite{article:Harten83,harten1984class}, Sweby \cite{sweby1984high},
and Tadmor \cite{tadmor1988convenient}
who developed mathematical precise limiting strategies
based on the principle that the numerical scheme should be 
total variation non-increasing. A good review of these limiters can be found in
Chapter 6 of LeVeque \cite{book:Le02}. 
In more recent work, many contributions have been made to
 generalize these limiting strategies to high-order methods, including to
 weighted non-oscillatory (WENO) schemes 
 (e.g., see Liu et al. \cite{LiuOsherChan94} and Jiang and Shu \cite{JiangShu96} ) and discontinuous Galerkin schemes 
 (e.g., for a review see Qiu and Shu \cite{article:QiuShu05} and Dumbser et al. \cite{dumbser2014posteriori}).

In this work we focus on limiters for discontinuous Galerkin (DG) methods \cite{ReedHill73,CoShu98,CoShu01}. 
Limiting strategies for DG methods can generally be broken down into two key steps:
\begin{description}
\item[{\bf Step 1.}] Detect {\it troubled-cells} (i.e., cells that contain large
gradients);
\item[{\bf Step 2.}]  Apply a limiter to each troubled cell that is detected
(i.e., in some way reduce the local gradients or higher derivative analogs thereof without changing the
average value in that cell). 
\end{description}
In some limiting strategies these two steps are executed with completely 
different techniques (e.g., see \cite{article:QiuShu05}), while in other approaches
the detection and limiting happen in all in one fell swoop (e.g., the moment limiter
of Krivodonova \cite{Kri07}).
In either case, the end goal of most limiters is to design a method that simultaneously accomplishes two
often conflicting objectives: (1) reduce unphysical oscillations in the presence of shocks and 
(2) retain high-order accuracy  in smooth regimes. 

The techniques used for DG methods in the limiting step ({\bf Step 2}) can be subdivided
into three broad classes of limiters, which we briefly describe below.

\subsubsection {Moment or slope methods}  These limiters have their roots 
in the development of classical second-order finite volume methods.
The key ingredients are minimum modulus evaluations that compare
the slope (or higher derivative analogs) in the current cell to values computed from
neighboring cells. Such methods are still widely used (e.g., see 
\cite{HotAckMosPhi04,hoteit2004new,Kri07,YangWang09,kuzmin2010vertex,kuzmin2013slope,kuzmin2014hierarchical}).
There are two well-known difficulties with these limiters: (1) 
they tend to degrade the of order of accuracy in smooth extrema; and (2)
they are difficult to generalize to unstructured meshes.

\subsubsection {Hybrid weighted essentially non-oscillatory methods}  
Methods in this category
are motivated by the success of the weighted essentially non-oscillatory (WENO)
method in shock-capturing.  The classical schemes
\cite{QiuShu04,QiuShu05,JQiuShu05,ZhuQiuShuDu08,ZhuQiu09,ZhuQiu13,ZhuZhoShuQiu13} redefine the
polynomial representation of the solution inside each element by considering
neighboring elements and matching cell averages.  The WENO
reconstruction, which is automatically mass conservative, is then typically applied in order 
to redefine polynomial representations inside each element.  More recently,
the so-called Hermite WENO schemes have been used in an effort to
reduce the size of the computational stencil
\cite{LuoBaumLohner07}.
A chief criticism of these methods has been the expense of their application;
although recent efforts exist to reduce their computational cost \cite{ZhuZhoShuQiu13}.

\subsubsection {Artificial viscosity methods} This method is based on
the classical idea of adding a small amount of artificial viscosity to
hyperbolic equations \cite{article:NeumRicht50}. The modern
incarnation of this limiting strategy was pioneered by Persson and Peraire \cite{Persson06}
and further developed in \cite{YuYan13,AbbMashJac14}.
These methods typically work by using a smoothness indicator to first detect troubled cells,
and then directly discretizing a diffusive term that has been artificially added to the equations in an effort to smooth-out oscillations.  One drawback of using these
limiters with explicit time-stepping schemes is that due to the diffusion operator, the
maximum stable time-step is $\dt = \BigOh( \hgrid/\degtwo )$
(i.e., the viscosity parameter is typically $\BigOh(\hgrid)$), 
%
%
where $\hgrid$ is the
grid spacing and $\deg$ is the highest polynomial order of the basis polynomials,
rather than the typical $\dt = \BigOh( \hgrid/\deg )$ for hyperbolic systems.
%

Other variations such as minimum entropy satisfying limiters
\cite{article:ZhangShu12} have also begun to be investigated.
A complete description of what options exist for limiting DG solutions is beyond the scope of this work; we direct the interested reader
to \cite{dumbser2014posteriori} for a recent and extensive review of the literature.

\subsection{Brief summary of current work}
The purpose of this work is to develop a novel limiter for discontinuous Galerkin 
methods that has all of the following properties:
\begin{itemize}
	\item Provable high-order accuracy in smooth regimes;
	\item Robust reduction in spurious oscillations in the presence of shocks;
	\item Computationally efficient (i.e., only uses neighboring cells and no
		  characteristic decompositions);
	\item Easy to implement in existing DG codes;
	\item Automatic extension to multidimensional settings including Cartesian and unstructured grids; and
	\item Automatic incorporation of positivity preserving limiters.
\end{itemize}
The limiter presented in this work can be viewed as a novel extension of two separate methods.
First, it can be thought of as extending the finite volume
Barth-Jespersen limiter \cite{BarthJesp89} to the discontinuous Galerkin framework.
Second, it can be viewed as an extension of the modern maximum principle preserving DG
schemes developed by Zhang and Shu \cite{ZhangShu11}.
In particular, the limiter is based on finding for each element local upper and lower bounds using only nearest neighbor values, and then limiting the higher-order modes
in such a way that that mass conservation is automatic and
such that the limited solution remains within the predicted local bounds.
Roughly speaking, this limiter can be viewed as belonging to the {\bf moment
or slope methods}, although it has important differences with the current methods
in this class, including no need to for characteristic decompositions, no need
for hierarchical highest-to-lowest moment limiting, and easy extensibility to
unstructured meshes.

The remainder of this paper is organized as follows.  In Section \ref{sec:dg-fem-framework}
 we briefly review the DG method. The basic limiting strategy for a one-dimensional
 scalar equation is outlined in Section \ref{sec:limitintro}.
 In Section \ref{sec:importance_of_alpha} we describe how the proposed limiter obtains
high-order accuracy. Extension to  multiple dimensions is shown in
Section \ref{sec:importance_of_phi} and to general systems of hyperbolic
 conservation laws in Section \ref{sec:systems}, where we use 
 the compressible Euler equations as the prototypical example.
In Section \ref{sec:positivity-preservation} we describe how
positivity-preserving ideas can be incorporated. 
The numerical method is validated on several standard cases for the
compressible Euler equations in Section \ref{sec:numerical-results}.
Conclusions are made in Section \ref{sec:conclusions}.

\subsection{Implementation details}
All of the numerical examples carried out in this work are implemented in the
software package {\sc DoGPack} \cite{dogpack}. This software package is
open-source and all the code used in this work can be freely downloaded from
the web.

\section{DG-FEM framework}
\label{sec:dg-fem-framework}

In order to define the proposed limiter, we first 
provide a brief description of the discontinuous
Galerkin (DG) method, which serves to define the notation used throughout the
remainder of this work.
We refer the interested reader to \cite{CoShu01,HesWar07,seal12} and
references therein for further details of the DG method.

Consider a conservation law of the form
\begin{equation}
\label{eqn:conslaw}
    q_{,t} + {\nabla} \cdot {\bf F}(q) = 0,
\quad \text{in} \, \, \, \Omega \subset \reals^d,
\end{equation}
with appropriate boundary conditions,
where $q(t,{\bf x}): \reals_{\ge 0} \times \reals^d \mapsto \reals^{\meq}$ are
the conserved variables, 
${\bf F}(q): \reals^{\meq} \mapsto \reals^{\meq \times d}$ is the flux function,
$d$ is the spatial dimension, and $\meq$ is the number of equations
in the system.
We assume that the system is {\it hyperbolic}, meaning that the flux Jacobian matrix,
${\bf n} \cdot {\bf F}_{,q}$, is diagonalizable with real eigenvalues for
all $q$ in the domain of interest and
for all $\| {\bf n} \| = 1$.

Let $\Omega \subset \reals^d$ be a polygonal domain with boundary $\partial \Omega$. 
The domain $\Omega$ is discretized via a finite set of non-overlapping elements, $\Tm_i$, such
that  $\Omega = \cup_{i=1}^N \Tm_i$. 
Let ${P}^{\, \deg}\left(\reals^d\right)$ denote the set of polynomials from $\reals^d$ to $\reals$ 
with maximal polynomial degree $\deg$.
Let $\WS^{h}$ denote the {\it broken} finite element space on the mesh:
\begin{equation}
\label{eqn:broken_space}
   \WS^h := \left\{ w^h \in \left[ L^{\infty}(\Omega) \right]^{\meq}: \,
    w^h \bigl|_{\Tm_i} \in \left[ {P}^{\, \deg} \right]^{\meq}, \, \forall \Tm_i \in \Tm^h \right\},
\end{equation}
where $h$ is the grid spacing.
The above expression means that $w^h \in \WS^{h}$ has
$\meq$ components, each of which when restricted to some element $\Tm_i$
is a polynomial of degree at most $\deg$ and no continuity is assumed
across element edges (or faces in 3D). 

The approximate solution on each element $\Tm_i$ at time $t=t^n$ is of the form
\beq
q^h(t^n, {\bf x}(\bxi)) \Bigl|_{\Tm_i} = \sum_{\ell=1}^{\mleg(\deg)} \, Q^{(\ell) n}_i  \,
\varphi^{(\ell)}\left( \boldsymbol{\xi} \right),
\eeq
where $\mleg$ is the number of Legendre polynomials and
$\varphi^{(\ell)}\left( \bxi \right): \reals^d \mapsto \reals$ are the Legendre polynomials defined
on the reference element ${\mathcal T}_0$ in terms of the reference coordinates $\bxi \in \Tm_0$.
The Legendre polynomials are orthonormal with respect to the following inner product:
\beq
\frac{1}{|\Tm_0|} \int_{\Tm_0} \varphi^{(k)}(\bxi) \, \varphi^{(\ell)}(\bxi) \, d\bxi = 
\begin{cases}
1 & \, \text{if} \quad k=\ell, \\
0 & \, \text{if} \quad k \ne \ell.
\end{cases}
\eeq

We note that independent of $h$, $d$, $\deg$,  and the type of element, the lowest order Legendre polynomial is always $\varphi^{(1)} \equiv 1$. This makes the first Legendre coefficient the
cell average:
\beq
    Q^{(1) n}_i = \frac{1}{|\Tm_0|} \int_{\Tm_0}
    q^h(t^n, {\bf x}(\bxi)) \Bigl|_{\Tm_i}  \, \varphi^{(1)}\left( \bxi \right) \,  d\bxi =: \bar{q}^n_i.
\eeq

To obtain the semi-discrete DG method, equation \eqref{eqn:conslaw} is
multiplied by the test function $\varphi^{(k)}$, the resulting equation
is integrated over an element $\Tm_i$,
and integrations-by-part are performed on the flux term. The result is
the following system of ODEs:
\begin{equation}
\label{eqn:semidiscrete_dg}
\dot{Q}^{(k)}_i = \frac{1}{|\Tm_i|} \int_{\Tm_i} \nabla_{\bf x} \varphi^{(k)} \cdot
 {\bf F}\left( q^h\Bigl|_{\Tm_i} \right) \, d{\bf x}
-  \frac{1}{|\Tm_i|} \oint_{\partial_{\Tm_i}}  \varphi^{(k)} {\mathcal F} \cdot d{\bf s},
\end{equation}
where the numerical fluxes ${\mathcal F}$ must be determined via 
(approximate) Riemann solvers. In this work we make use of the Rusanov
 \cite{article:Ru61} numerical flux, which is often referred to as the  local Lax-Friedrichs (LLF) flux in the literature.

Finally, in order to obtain a Runge-Kutta DG method,
equation \eqref{eqn:semidiscrete_dg} is discretized in time via some Runge-Kutta method.
In this work we make use of strong stability-preserving (SSP) Runge-Kutta methods,
such as those found in Gottlieb, Shu, and Tadmor \cite{gottliebShuTadmor01} and Ketcheson \cite{Ke08}.

\section{Proposed limiter: The scalar case}
\label{sec:limitintro}
The basic procedure we follow to limit the solution on a single element $\Tm_i$ 
is given by the following sequence of steps:

\begin{description}
 \item[{\bf Step 0.}] Define a set of points $\chi_i$ that will be used to approximate cell
 maximum and minimum values.  
 In the current work, we use Gaussian quadrature nodes and  augment them with corner points and 
 quadrature points along the element boundaries (i.e., edge Gaussian quadrature points).
%
\item[{\bf Step 1.}] For each mesh element, $\Tm_i$, we compute an approximate maximum and minimum:
\begin{equation*}
    q_{M_i} := \max_{x\in\chi_i} \left\{ q^{h}( x )\Bigl|_{\Tm_i} \right\}
\quad \text{and} \quad
    q_{m_i} := \min_{x\in\chi_i} \left\{ q^{h}( x )\Bigl|_{\Tm_i} \right\}.
\end{equation*}
\item[{\bf Step 2.}] We consider the set $N_{\Tm_i}$ of all neighbors of $\Tm_i$, 
 excluding $\Tm_i$ itself, and compute an approximate upper and lower bound:
 \begin{align}
     M_i :=& \max\left\{\bar{q}_i+\alpha(\hgrid), \, \, \max_{j \in N_{\Tm_i}}\left\{q_{M_j}
     \right\} \right\}, \\
     m_i :=& \min\left\{\bar{q}_i-\alpha(\hgrid), \, \, \min_{j \in N_{\Tm_i}}\left\{q_{m_j}
     \right\}
     \right\}.
 \end{align}
 The scalar function $\alpha(\hgrid)\geq0$ is a tolerance function that will be
 described shortly.  
 The most aggressive limiter we consider sets 
 $\alpha = 0$.
\item[{\bf Step 3.}] We define 
 \begin{equation}
 \label{eqn:theta}
 	\theta_{M_i} := \phi\left( \frac{ M_i-\bar{q}_i }{ q_{M_i}-\bar{q}_i }\right)
	\quad \text{and} \quad
    \theta_{m_i} := \phi\left( \frac{ m_i-\bar{q}_i }{ q_{m_i}-\bar{q}_i }\right),
 \end{equation}
 where $0 \leq \phi(y) \leq 1$ is a cutoff function.  In this work, we
 use the function 
 \begin{equation}
 \phi(y) := \min\left\{ \frac{y}{1.1}, 1 \right\},
 \end{equation}
which is motivated by the work of \cite{michalak2008limiters}.
The form of this function becomes important in the multidimensional setting, which will be described shortly.
\item[{\bf Step 4.}] We define the \emph{rescaling parameter} as
 \begin{equation}
 \label{eqn:thetai}
     \theta_{i}:=\min\left\{ 1,\ \theta_{m_i},\ \theta_{M_i} \right\}.
 \end{equation}
 \item[{\bf Step 5.}] Finally, we rescale the approximate solution on the element $\Tm_i$ as
 \begin{equation}
 \label{eqn:rescaling}
    \tilde{q}^h( x )\Bigl|_{\Tm_i} := \bar{q}_i + \theta_i \left( q^h(x)\Bigl|_{\Tm_i} - \bar{q}_i  \right).
 \end{equation}
 
\end{description}


Before delving into the finer details of the method, a few remarks are in order.

\begin{rmk}
The choice of quadrature points in Step 1 is not unique.
\end{rmk}
Indeed, the purpose of this step is to construct an approximate upper and lower bound 
for the solution.  We originally tested this limiter by using cell averages,
but found that those solutions tended to be quite diffusive near shocks.
Directly sampling neighboring cells (that are potentially oscillating) allows us to retain 
sharper features in the solution.

\begin{rmk}
The presence of a non-zero value of $\alpha$ in Step 3 is required to obtain high-order accuracy
at local extrema.  
\end{rmk}
This issue will be elaborated upon in \mbox{Section
\ref{sec:importance_of_alpha}}.

\begin{rmk}
When the proposed scheme is compared to the recent maximum principle preserving (MPP) limiters \cite{ZhangShu11},
significant differences come into play in {\bf Step 2}  and {\bf Step 3}.
\end{rmk}
\begin{itemize}
\item For the MPP methods, {\bf Step 2} is replaced by
fixed {\it a priori} {global} bounds, whereas we follow the 
Barth-Jespersen idea \cite{BarthJesp89} and use this step to estimate
\emph{local} upper and lower bounds for the solution.
\item In {\bf Step 3}, the MPP limiters use hard upper and lower bounds
for the solution that are known \emph{a priori}, whereas we use the result
from {\bf Step 2} to estimate these.
\item In the MPP literature, authors normally set $\phi( y ) = \min\{ 1, y \}$, whereas
we find that the form of this function is important when pushing to multiple dimensions.
In Section \ref{sec:importance_of_phi} we elaborate on this idea.
\end{itemize}
Additionally, the newly proposed limiter is designed to  capture shocks,
whereas the original MPP limiter for DG methods was solely designed to
preserve the maximum principle property
of hyperbolic problems.  In their work, additional
limiting was necessary in order to suppress unphysical oscillations.

\begin{rmk}
The rescaling in Step 6 does not effect the total mass.
\end{rmk}

The consequence of this is that the limiter will be automatically mass conservative, which
is an important property for hyperbolic solvers.  This can be observed by simply integrating
equation \eqref{eqn:rescaling} over a single element.

\section{On the importance of $\alpha$ for retaining high-order accuracy}
\label{sec:importance_of_alpha}

In this section we describe the importance of $\alpha$ for obtaining a method
that is genuinely high-order accurate.  Many limiters exhibit the so-called 
``clipping" phenomenon, where smooth extrema are cut off because the limiter
turns on and damps out these values.  
Recent efforts have been put forth to mitigate this clipping
phenomenon for other discretizations, including
the piecewise parabolic method (PPM) \cite{ColSek08}.

In this section, we will argue that high-order accuracy will 
be retained provide that the scalar function $\alpha\left(\hgrid \right) : \R_{\geq 0} \to \R_{\geq 0}$ vanishes
slower than $\BigOh( \hgrid^2 )$.  This observation was inspired by recent work on
extensions of the Barth-Jespersen limiter to finite volume methods \cite{michalak2008limiters},
and we demonstrate the performance of our limiter by testing it on a smooth solution in 1D.

In one dimension we can break down any smooth function on a bounded domain 
into monotonic regions and finitely many points of local extrema.  We begin with
a discussion of monotonic regions, and then turn to isolated local extrema.
\subsection{Performance in monotonic regions}
\label{subsec:monotonic}
In regions where the exact solution is monotonic, the proposed limiter
does not effect the asymptotic convergence rate independent of the choice of $\alpha \geq 0$.
Without loss of generality, we consider a region with the solution is non-decreasing.
After after enough refinement, the exact solution satisfies
$q_{i-1} \le q_i \le q_{i+1}$ 
where $q_{i}$ refers to the solution
value on the $i^{\text{th}}$-grid cell in a monotonic region.  The algorithm then correctly
selects $q_{i+1}$ as an upper bound, $q_{i-1}$ as a lower bound, and therefore does not overreach by over limiting
the solution on the $i^{\text{th}}$ grid cell.  A formal proof of the order of accuracy
would follow that already presented in \cite{ZhangShu11}, but would start with the fact that
our algorithm over-estimates the upper  bound and under-estimates
the lower bound for the exact solution.

%
%

\subsection{Performance in regions near smooth extrema}
In this subsection, we describe the importance of using a non-zero value of
$\alpha$ to maintain genuine high-order accuracy.  If $\alpha$ is chosen to be too
large, then oscillations in a given cell will never be clamped down, whereas if $\alpha$
vanishes too quickly, then smooth extrema will be clipped.  Our goal is to find
necessary conditions for maintaining high-order accuracy at smooth extrema.


To begin, we consider a Taylor expansion of a smooth function $q : \R \to \R$ at an 
extrema $x=\xi_0$:
\begin{equation*}
    q(x)=q({\xi_0}) + q'(\xi_0) (x-\xi_0)  + \frac{1}{2} q''(\xi_0) (x-\xi_0)^2 + \cdots.
\end{equation*}
Given that $q'(\xi_0) = 0$, we observe that asymptotically $\left|q(x)- q( {\xi_0} ) \right| = \BigOh(\hgrid^2)$
if  $\left|x- {\xi_0}  \right| = \BigOh(\hgrid)$.
This leads us to consider functions of the form 
$\alpha (\hgrid) = \BigOh(\hgrid^r)$, where $r < 2$ in order to maintain
\begin{equation*}
    |q(x)-q( {\xi_0} )|  = \BigOh( \hgrid^2 )\le |\alpha(\hgrid)| \quad \text{for all} \quad | x - \xi_0 | = \BigOh(\hgrid).
\end{equation*}
Note that in the 2D Cartesian case: $h = \text{max} \left( \dx, \, \dy \right)$,
where $\dx$ and $\dy$ are the mesh widths in each coordinate direction.


That is, as long as $\alpha(\hgrid)$ goes to zero \emph{slow enough}, then our limiter
will eventually turn off and not clip smooth extrema.  However, if we were to set
$\alpha = 0$, then our limiter will induce clipping at smooth extrema, and yet if $\alpha$
does not vanish as $h \rightarrow 0$, then our method would permit problematic oscillations to persist.

The choice of $\alpha$ is not unique, and if desired, could be tuned for any given problem.  The trade-off
is that if $\alpha$ is too large, then true oscillations will only be damped
provided the mesh is fine enough.
Our choice of $\alpha(\hgrid) = \BigOh( \hgrid^{1.5} )$ is consistent with the recommended switches that turn off similar Barth-Jespersen
type limiters for finite volume schemes \cite{michalak2008limiters}.  Extensive numerical tests 
indicate this choice finds a good balance between the need to limit near shocks,
while being soft enough to allow high-order accuracy on coarse grids.  
{\bf In the case of smooth problems, our limited solution will always converge to the unlimited
solution even near local extrema.}
We illustrate this point with a simple case study.

\subsection{On the importance of $\alpha$: A simple 1D case study}
\label{subsec:1dadvection}
We analyze the performance of the proposed limiter on the
linear advection equation \eqref{eqn:advection1d} on $x \in [0,1]$ with
periodic boundary conditions and with the initial condition
\begin{equation}
    q(x) = \begin{cases}
        \cos^6\left( \frac{ (x-0.5) \pi}{0.16} \right) \quad & \text{if } |x-0.5| < 0.08,
    \\ 0 \quad & \text{otherwise}.
    \end{cases}
\end{equation}
In Table \ref{table:alpha_convergence} we present a convergence study for various
grid resolutions after advecting the solution through one full period (e.g., solve
to $t=1$).  
For this example we use the $10$-stage SSP fourth-order Runge-Kutta scheme (SSPRK4)
developed by Ketcheson \cite{Ke08} and hold a constant CFL number of $0.4$.
We compare the results for the unlimited solution,
the solution with $\alpha = 0$, and the solution with $\alpha(\hgrid) = 50\hgrid^{1.5}$.

We observe that after enough refinement, the limited and unlimited solutions
agree with each other, whereas for very coarse grids, the order of accuracy is slightly
degraded.  It is only after the limiter ``turns off'' that we
we dial in on high-order accuracy. 
This cutoff point is a function of the curvature of the exact solution at the
local extrema.
This example serves to illustrate that after enough refinement, there is no
difference between the limited and unlimited solutions.

 \begin{table}[ht]
 \begin{center}
\hspace*{-50pt}\begin{tabular}{|r||c|c||c|c||c|c|}
\hline
\bf{Mesh} & \bf{{No Limiter}} & \bf{{Order}}& \bf{{$\alpha=0$}} & \bf{{Order}}& \bf{{$\alpha=50 \hgrid^{1.5}$}} & \bf{{Order}}\\
\hline
\hline
$   5$ & $3.78\times 10^{-01}$ & --- & $8.73\times 10^{-01}$ & --- & $3.78\times 10^{-01}$ & ---\\
\hline
$   8$ & $5.61\times 10^{-01}$ & $-0.84$ & $9.00\times 10^{-01}$ & $-0.07$ & $5.61\times 10^{-01}$ & $-0.84$ \\
\hline
$  14$ & $2.32\times 10^{-01}$ & $1.58$ & $8.22\times 10^{-01}$ & $0.16$ & $2.34\times 10^{-01}$ & $1.56$ \\
\hline
$  24$ & $5.38\times 10^{-02}$ & $2.71$ & $7.04\times 10^{-01}$ & $0.29$ & $1.74\times 10^{-01}$ & $0.56$ \\
\hline
$  42$ & $4.24\times 10^{-03}$ & $4.54$ & $4.78\times 10^{-01}$ & $0.69$ & $7.82\times 10^{-02}$ & $1.43$ \\
\hline
$  73$ & $3.65\times 10^{-04}$ & $4.44$ & $2.24\times 10^{-01}$ & $1.38$ & $2.19\times 10^{-02}$ & $2.30$ \\
\hline
$ 127$ & $3.89\times 10^{-05}$ & $4.04$ & $8.84\times 10^{-02}$ & $1.68$ & $2.73\times 10^{-03}$ & $3.76$ \\
\hline
$ 222$ & $4.10\times 10^{-06}$ & $4.03$ & $2.97\times 10^{-02}$ & $1.95$ & $4.10\times 10^{-06}$ & $11.64$ \\
\hline
$ 388$ & $4.38\times 10^{-07}$ & $4.01$ & $9.37\times 10^{-03}$ & $2.07$ & $4.38\times 10^{-07}$ & $4.01$ \\
\hline
$ 679$ & $4.65 \times 10^{-08}$ 
& $4.01$ & $2.76\times 10^{-03}$ & $2.19$ & $4.65\times 10^{-08}$ & $4.01$ \\
\hline
\end{tabular}\hspace*{-50pt}

\caption{ 
    Convergence for the 1D advection problem in Section \ref{subsec:1dadvection}. All errors are $L_2$ norm errors. We
    see that a nonzero value of $\alpha$ is required to allow the limiting to completely turn off when the solution
    is properly resolved.}
\label{table:alpha_convergence}

\end{center}
\end{table}

\section{On the importance of $\phi$ for multidimensional cases}
\label{sec:importance_of_phi}

An important distinction between the single and multidimensional cases is that
in higher dimensions, shocks are no longer isolated to single cells.  
To illustrate this, consider the following initial condition:
\begin{gather}
	q(x,y) = \begin{cases}
	1, & \text{if } x \leq 0, \\
	0, & \text{otherwise}.
	\end{cases}
\end{gather}
Now, consider a Cartesian mesh with cells of width $\dx$ and height $\dy$:
\begin{equation}
	C_{ij} := \left[ \left(i-\half \right)\dx,  \left(i+\half \right)\dx \right] 
	\times \left[\left(j-\half\right)\dy, \left(j+\half\right)\dy \right], \quad i,j \in \Z.
\end{equation}
After projecting the initial conditions onto ${\mathcal W}^h$
with $\deg>1$ (i.e., see definition \eqref{eqn:broken_space}),  
every cell $C_{0j}$ $\forall j$ will 
exhibit Gibbs phenomena; this is because the initial 
discontinuity bisects the cells  $C_{0j}$ $\forall j$.
We now consider the effect of applying our limiter to this solution.
 We focus on what happens to the
solution on cell $C_{00}$.  

In {\bf Step 2} of the proposed limiter, we compute upper and lower bounds of the solution using nearest neighbors,
which in this case satisfy 
\[
	\max \left\{ q^h \bigl|_{C_{00}} \right\} = \max \left\{ q^h \bigl|_{C_{0 \, \pm 1}} \right\} \quad 
	\text{and} \quad
	\min q^h \left\{ \bigl|_{C_{00}} \right\} = \min \left\{  q^h \bigl|_{C_{0 \, \pm 1}} \right\}.
\]
That is, for these initial conditions, the predicted upper bound for cell $C_{00}$ will not clamp down
on any oscillations present in this cell, unless we are careful about deciding on the
correct structure of $\phi$.

The important observation to make here is that it is \emph{the ratio} of the deviation from cell averages
(the value $\theta$ in equation \eqref{eqn:theta} of {\bf Step 3}) that is the relevant
quantity to study.  For this problem, it is identically equal to 1 for every mesh element.
Indeed, what this tells us is that it is necessary for $\phi(1) < 1$ if we will have any hope
of clamping down on the unavoidable oscillations.  We observe that 
even with $\phi(1) < 1$, the oscillations are not 
completely suppressed after a single time step; however,
 after several time-steps the artificial oscillations decay exponentially to zero.
In practice, we observe that this is sufficient to construct robust simulations.
 
\begin{figure}[!t]
 \begin{center}
  \subfigure[$\phi(y) = \min\{1, y\}$]
    {\includegraphics[width=0.4\textwidth]{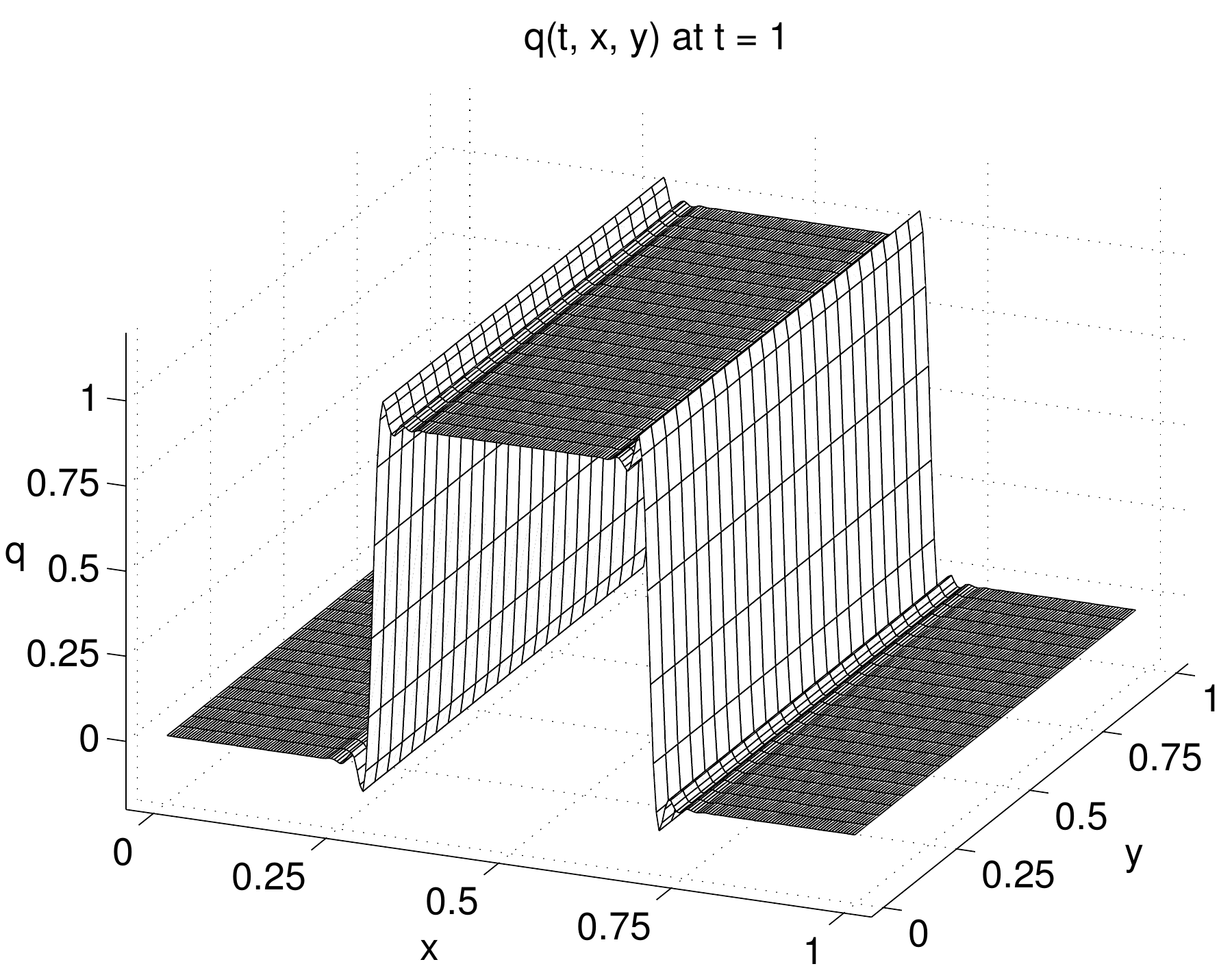}}
  \subfigure[1D slice of $\phi(y) = \min\{1, y\}$]
    {\includegraphics[width=0.4\textwidth]{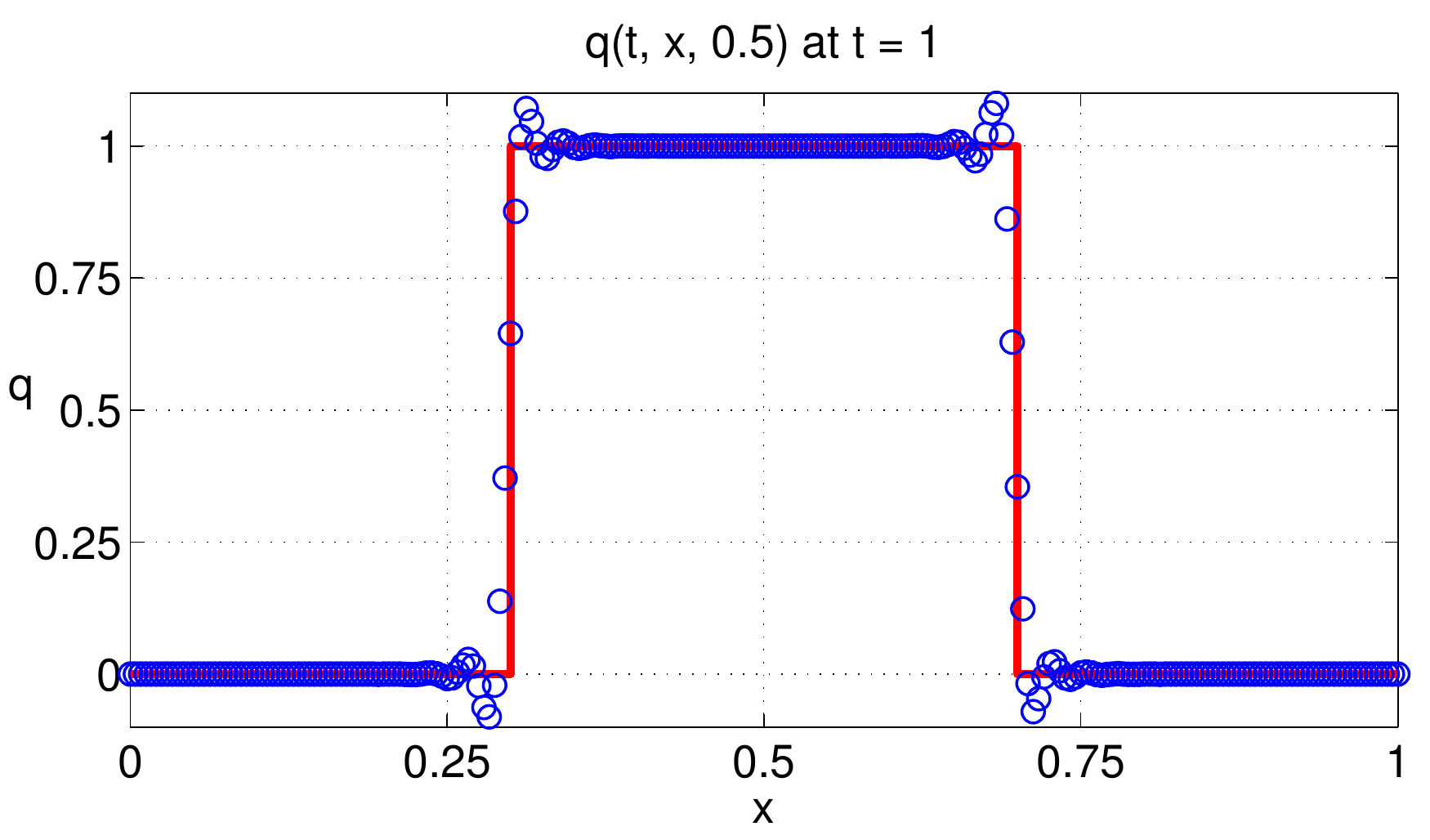}}
    \subfigure[$\phi(y) = \min\{1, \frac{y}{1.1} \}$]
    {\includegraphics[width=0.4\textwidth]{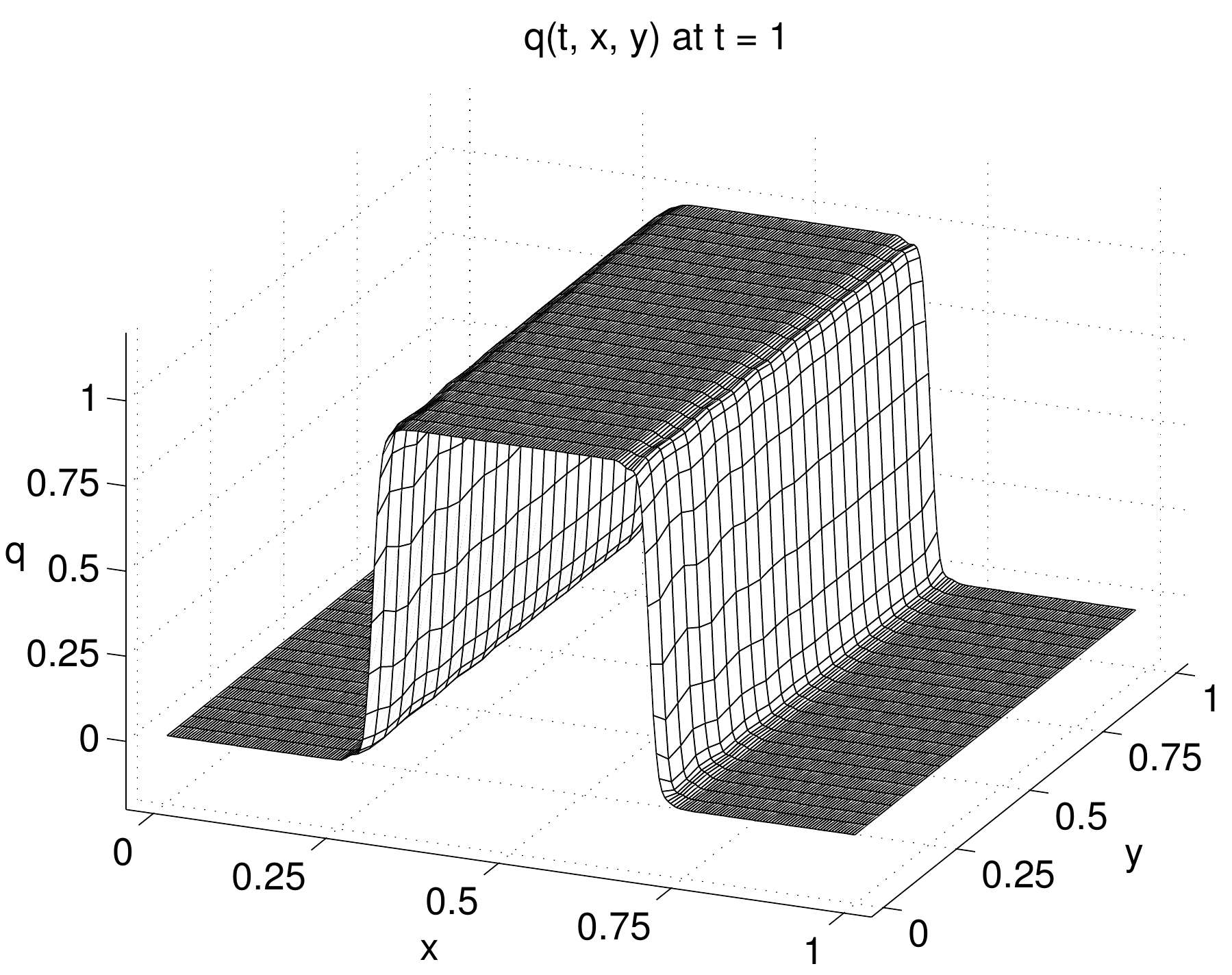}}
  \subfigure[1D slice of $\phi(y) = \min\{1, \frac{y}{1.1} \}$]
    {\includegraphics[width=0.4\textwidth]{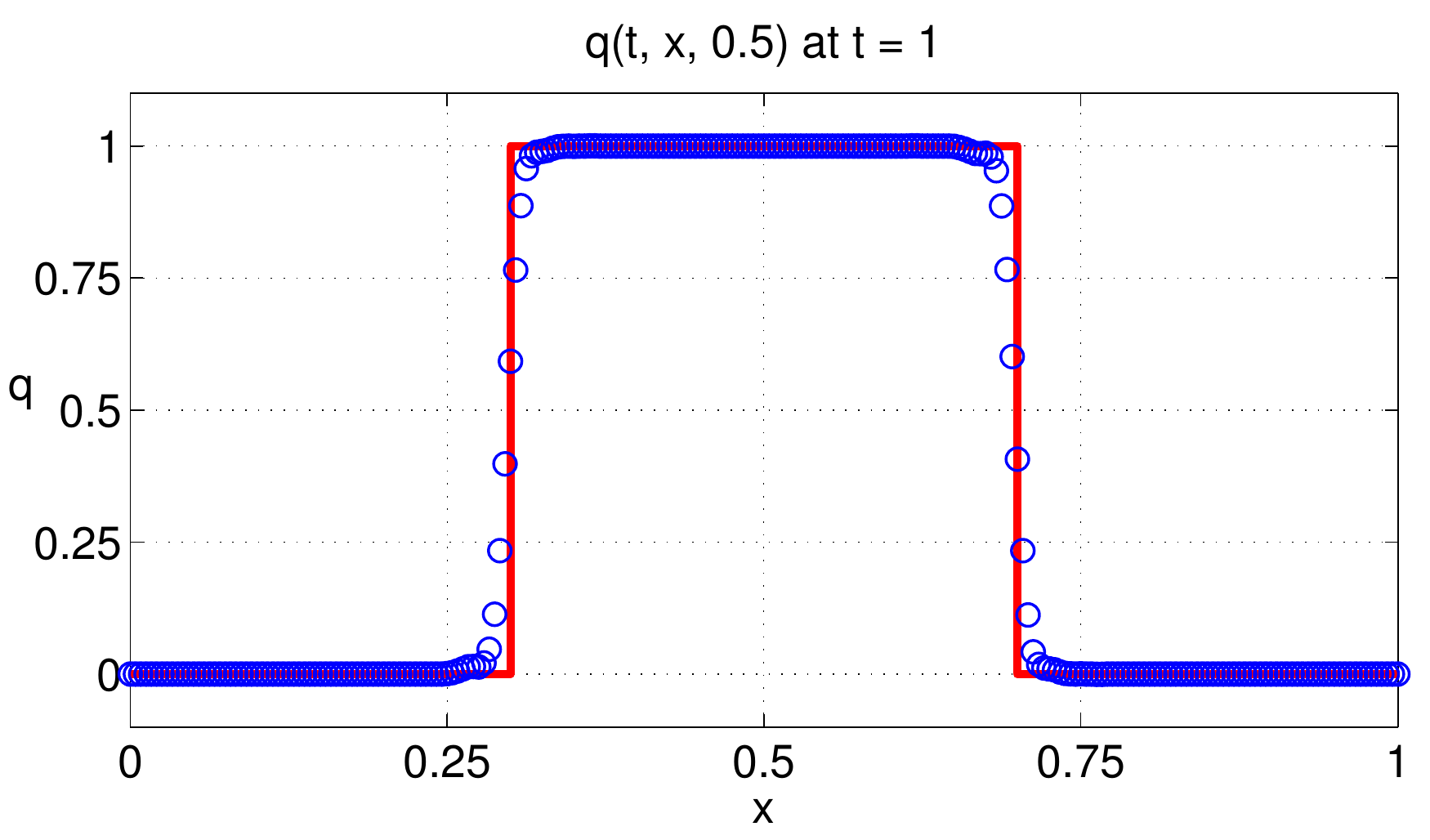}}
  \caption{Square wave test case for scalar advection with double periodic boundary
  conditions.  The solution in Panels
      (a) and (b) is computed using $\phi(y) = \min\{1, y\}$, while the solution
      in Panels (c) and (d) is computed using
      $\phi(y) = \min\{1, \frac{y}{1.1}\}$. 
       The simulation is run to
        a final time of $t=1$, at which point the solution should return to the
      initial conditions.  
      The choice in (a) and (b) satisfies $\phi(1) = 1$, whereas the choice in (c)
      and (d) satisfies $\phi(1) < 1$.  The first choice permits oscillations to remain in the
      solution, whereas the second choice eventually suppresses the unphysical oscillations.
      The remainder of the simulations used in this work will use the function
      $\phi(y) = \min\{1, \frac{y}{1.1}\}$.
  }
\label{fig:phi}
\end{center}
\end{figure}

As a simple demonstration of the importance of this choice of $\phi$, we consider linear advection 
in two-dimensions on a Cartesian grid.  Consider the scalar advection equation
\begin{equation}
	\frac{\partial q}{\partial t} + \frac{\partial q}{\partial x} + \frac{\partial q}{\partial y} = 0, \quad \text{for}  \quad (x,y) \in [0,1] \times [0,1],
\end{equation}
with double-periodic boundary conditions and discontinuous initial conditions:
\begin{equation}
    q_0(x,y) = \begin{cases}
    1 & \mbox{ if } x \in [0.3, 0.7], \\
    0 & \mbox{ otherwise}.
    \end{cases}
\end{equation}
In Figure \ref{fig:phi} we show results for a grid of size $60\times 60$
for a fourth-order ($\mathbb{P}_3$) method, comparing the choices of
 $\phi(y) = \min\{1, y\}$ and
$\phi(y) = \min\{ 1, \frac{y}{1.1} \}$. This result shows that
the second choice is necessary for pushing our results to multiple dimensions.
Other choices exist, but one property we find important is that $\phi(1) < 1$.
In the remainder of our numerical simulations, we choose $\phi(y) = \min\{ 1, \frac{y}{1.1} \}$
in all the 2D examples.

\section{Extensions to systems of equations}
\label{sec:systems}

As a prototypical example of a hyperbolic system of equations, we consider 
the compressible Euler equations:
\begin{gather}
\label{eqn:euler}
 \frac{\partial}{\partial t}
  \begin{bmatrix}
    \rho \\ \rho \ubold \\ \En
  \end{bmatrix} +  \nabla \cdot
  \begin{bmatrix} \rho \ubold \\ \rho \ubold  \ubold + {p} {\mathbb I} \\
    \ubold \left(\En + {p} \right) 
  \end{bmatrix} = 0,
  \qquad 
  \En = \frac{1}{2} \rho \| \ubold \|^2 + \frac{p}{\gamma-1},
\end{gather}
where $\gamma$ is the gas constant. The conserved and primitive variables are
\begin{equation}
    q := \left( \rho, \rho u^1, \rho u^2, \rho u^3, \mathcal{E} \right)
    \quad \text{and} \quad
    w := \left( \rho, u^1, u^2, u^3, p \right),
\end{equation}
respectively. In these expressions $\rho$ is the mass density,  
$\rho {\bf u} := (\rho u^1, \rho u^2, \rho u^3)$ is the momentum density, 
$\mathcal{E}$ is the total energy density, ${\bf u} := ( u^1,  u^2,  u^3)$  is the fluid velocity, and $p$ is the pressure.

The simplest implementation of the proposed limiter for  this system
 would be to apply the scalar limiting procedure to  each conserved variable independently.  Indeed, our initial tests employed this approach.
However, after extensive testing for this system,
we discovered that using
the primitive variables to define one high-order damping parameter $\theta$ yields better results. 
An example is provided in Section \ref{subsec:prim_vs_cons} to verify this claim.

The extension of the proposed limiter to systems of equations can
be summarized as follows: 
\begin{description}
\item[{\bf Step 0.}] Select a set of variables to use in the bounds checking: $w$.
In fluid dynamics, primitive variables are the preferred choice due to their
Galilean invariance. Select a set of points $\chi_i$ that will be used to approximate cell
 maximum and minimum values. In this work we always select: corners, internal, and edge Gaussian quadrature points.
\item[{\bf Step 1.}]  For each element $i$ and each component $\ell$ of $w$ compute:
\begin{equation}
    w^\ell_{M_i} := \max_{x\in\chi_i} \left\{ w^\ell(q^{h}( x ))\Bigl|_{\Tm_i} \right\}
     \quad \text{and} \quad
    w^\ell_{m_i} := \min_{x\in\chi_i} \left\{ w^\ell(q^{h}( x ))\Bigl|_{\Tm_i} \right\}.
\end{equation}
\item[{\bf Step 2.}] For each element $i$ and each component $\ell$ of $w$ compute an approximate upper and 
lower bound over the set of neighbors, $N_{\Tm_i}$
(excluding the current cell $\Tm_i$):
\begin{align}
\label{eqn:big_mi_sys}
     M^{\ell}_i &:= \max \left\{ \bar{w}^{\ell}_i + \alpha(\hgrid), \, \max_{j \in N_{\Tm_i}}
     \left\{ w^{\ell}_{M_j}\right\} \right\}, \\
\label{eqn:little_mi_sys}
     m^{\ell}_i &:= \min\left\{\bar{w}^{\ell}_i-\alpha(\hgrid), \,  \min_{j \in N_{\Tm_i}}
     \left\{w^{\ell}_{m_j}\right\} \right\},
 \end{align}
 where $\bar{w}^{\ell}_i$ is the cell average of ${w}^{\ell}$ over cell
 $\Tm_i$. $N_{\Tm_i}$.
 In the Cartesian grid case, we define $N_{\Tm_i}$ to be the set of cells that share a
 common edge with $\Tm_i$. This provides good results in general, and has the nice property that it does not extend
 the effective stencil size. 
 For unstructured grids, we define $N_{\Tm_i}$ as the set of all cells that
 share a common node with $\Tm_i$.  If this choice is replaced by elements that share common
 edges only, then the results are more diffusive.
 
\item[{\bf Step 3.}] For each element $i$ compute:
 \begin{equation}
 	\theta_{M_i} := \min_{\ell} \left\{ \phi\left( \frac{ M^{\ell}_i-\bar{w}^{\ell}_i }{ w^{\ell}_{M_i}-\bar{w}^{\ell}_i }\right) \right\}
	\, \, \, \text{and} \, \, \,
    \theta_{m_i} := \min_{\ell} \left\{ \phi\left( \frac{ m^{\ell}_i-\bar{w}^{\ell}_i }{ w^{\ell}_{m_i}-\bar{w}^{\ell}_i }\right) \right\}.
 \end{equation}
 \item[{\bf Step 4.}] For each element $i$ compute:
 \begin{equation}
     \theta_{i}:=\min\left\{ 1,\ \theta_{m_i},\ \theta_{M_i} \right\}.
 \end{equation}
 \item[{\bf Step 5.}] For each element $i$ and each component $\ell$ of $q$ compute 
 \begin{equation}
    \tilde{q}^h( x )\Bigl|_{\Tm_i} := \bar{q}_i + \theta_i \left( q^h(x)\Bigl|_{\Tm_i} - \bar{q}_i  \right).
 \end{equation}
\end{description}
Note that we are using the values of the variable $w$ to determine the
amount of limiting that needs to be done (i.e., size of $\theta_i$), but we are
actually applying the limiter to the conserved variables, thereby retaining mass conservation at the discrete level.

\subsection{Primitive vs. conservative variables for computing $\theta$: A 1D example}
\label{subsec:prim_vs_cons}
In order to test the limiter performance using either primitive or conservative variables to determine $\theta_i$, we consider the one-dimensional shock-entropy problem
found in \cite{ShuOsher89, LiuOsherChan94, SeGuCh14}, which
highlights the interplay between shock-capturing and feature-preservation.
The initial conditions for this problem are given by
\begin{align*}
   (\rho, u^1, p) =&
   \begin{cases} 
   \left( 3.857143, 2.629369, 10.3333 \right) & x < -4, \\
    \left( 1 + \varepsilon \sin(5x), 0, 1 \right) & x \geq -4.
   \end{cases}
\end{align*}
We follow the common practice of setting $\varepsilon = 0.2$ and discretize the
domain $[-5,5]$ into a total of 200 cells.  The final time for this simulation is $t=1.8$.
This test problem involves a shockwave interacting with a smooth entropy wave.
As this interaction occurs, the result is the formation of highly oscillatory
waves after the shock, where high-order accuracy should produce a benefit, yet limiters often smear out these oscillations. 
Ideally a limiter should be able to pick out and limit only non-physical
oscillations, while simultaneously capturing the oscillatory post shock
features.

\begin{figure}[!ht]
 \subfigure[]{\includegraphics[width=0.45 \linewidth]{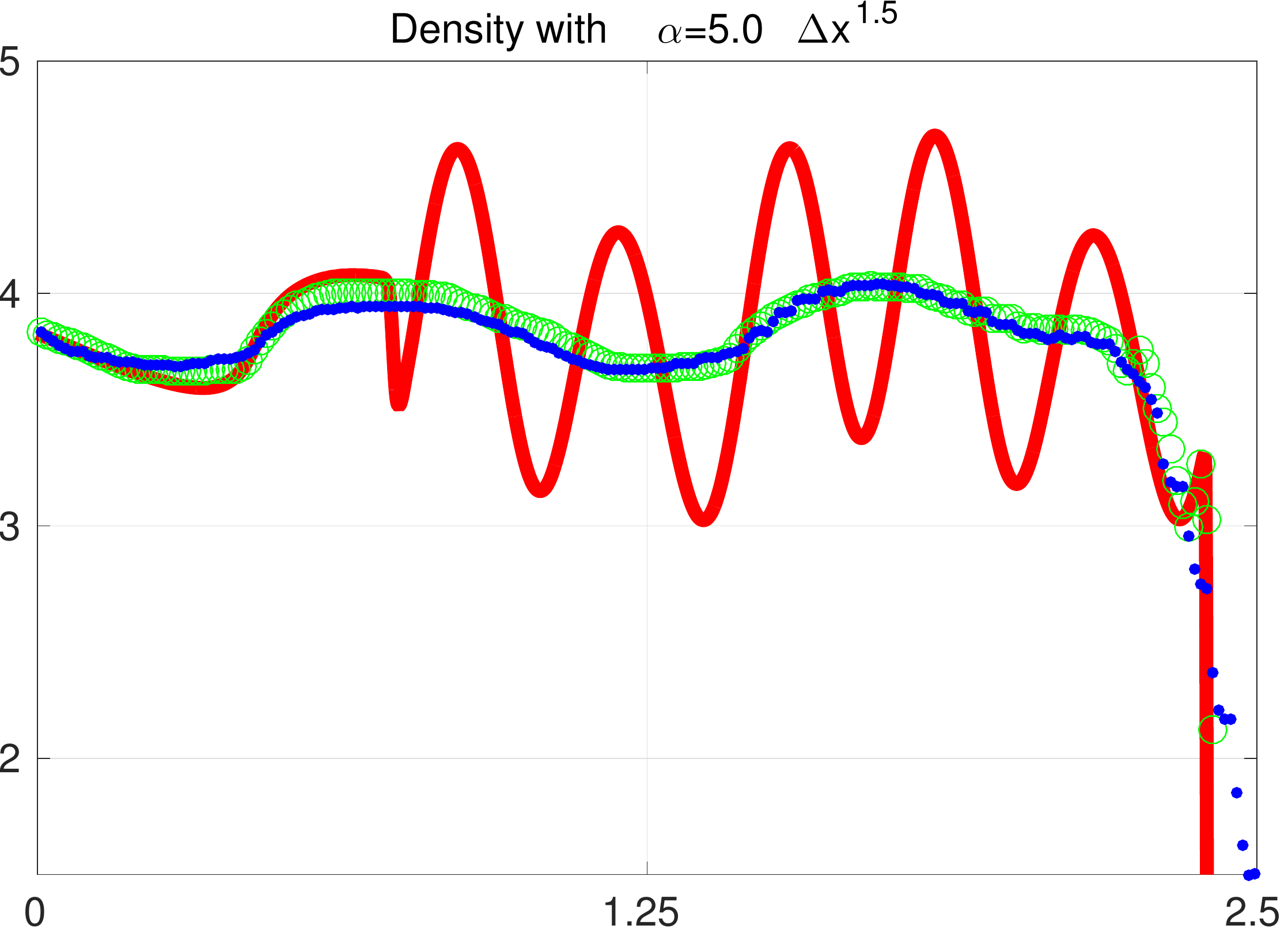}}
  \subfigure[]{\includegraphics[width=0.45 \linewidth]{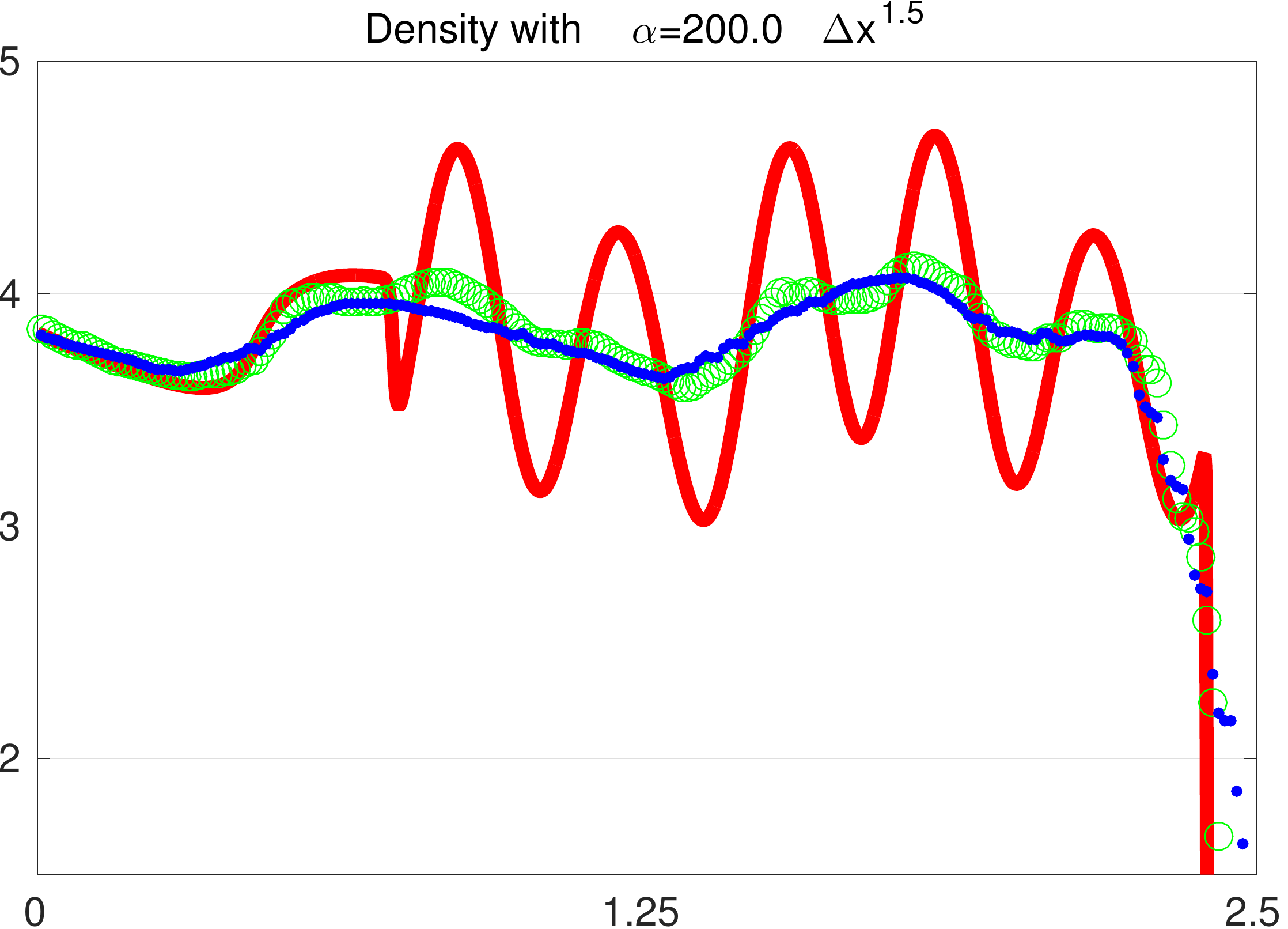}}
   \subfigure[]{\includegraphics[width=0.45 \linewidth]{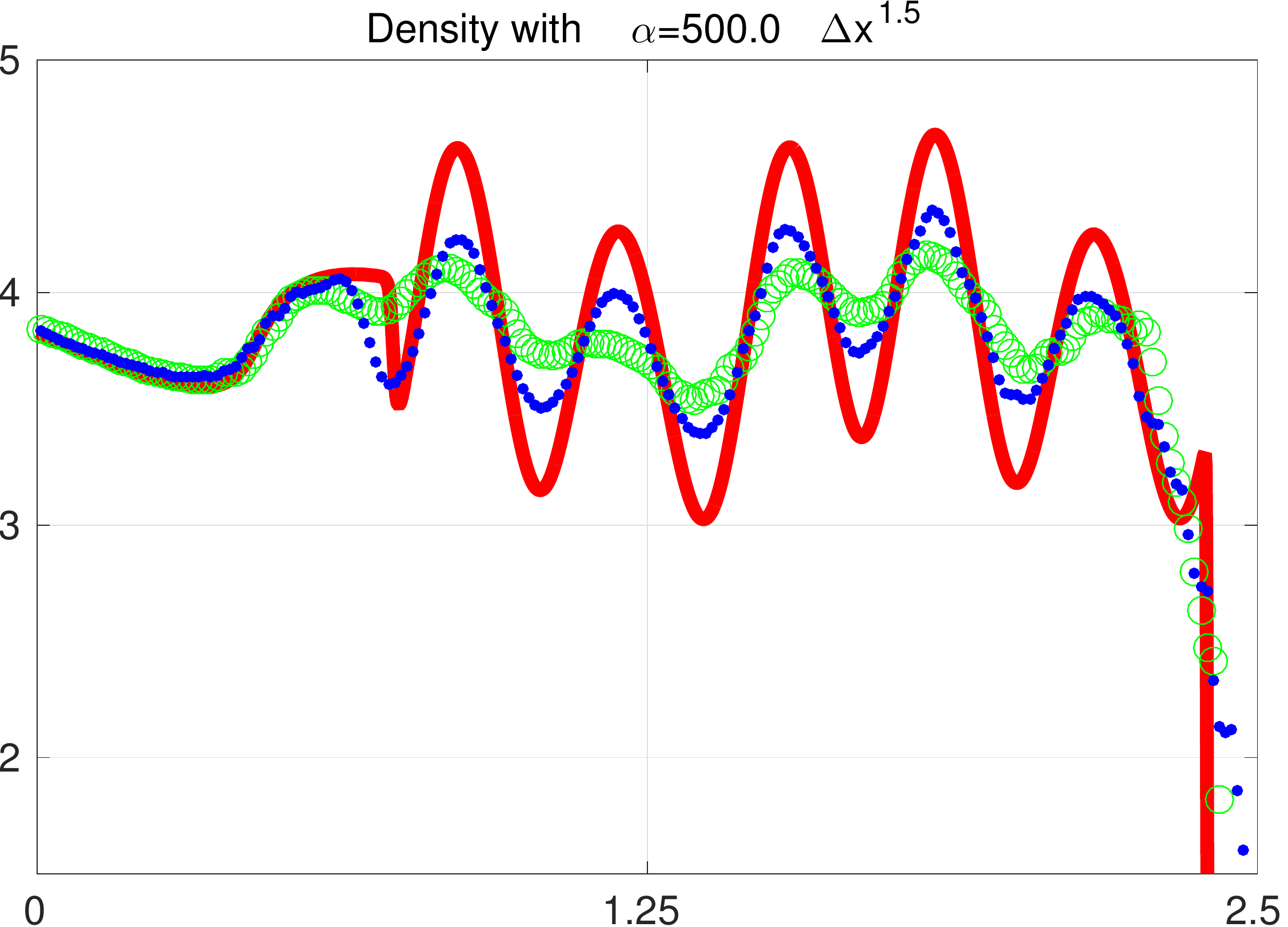}}
 \subfigure[]{\includegraphics[width=0.45 \linewidth]{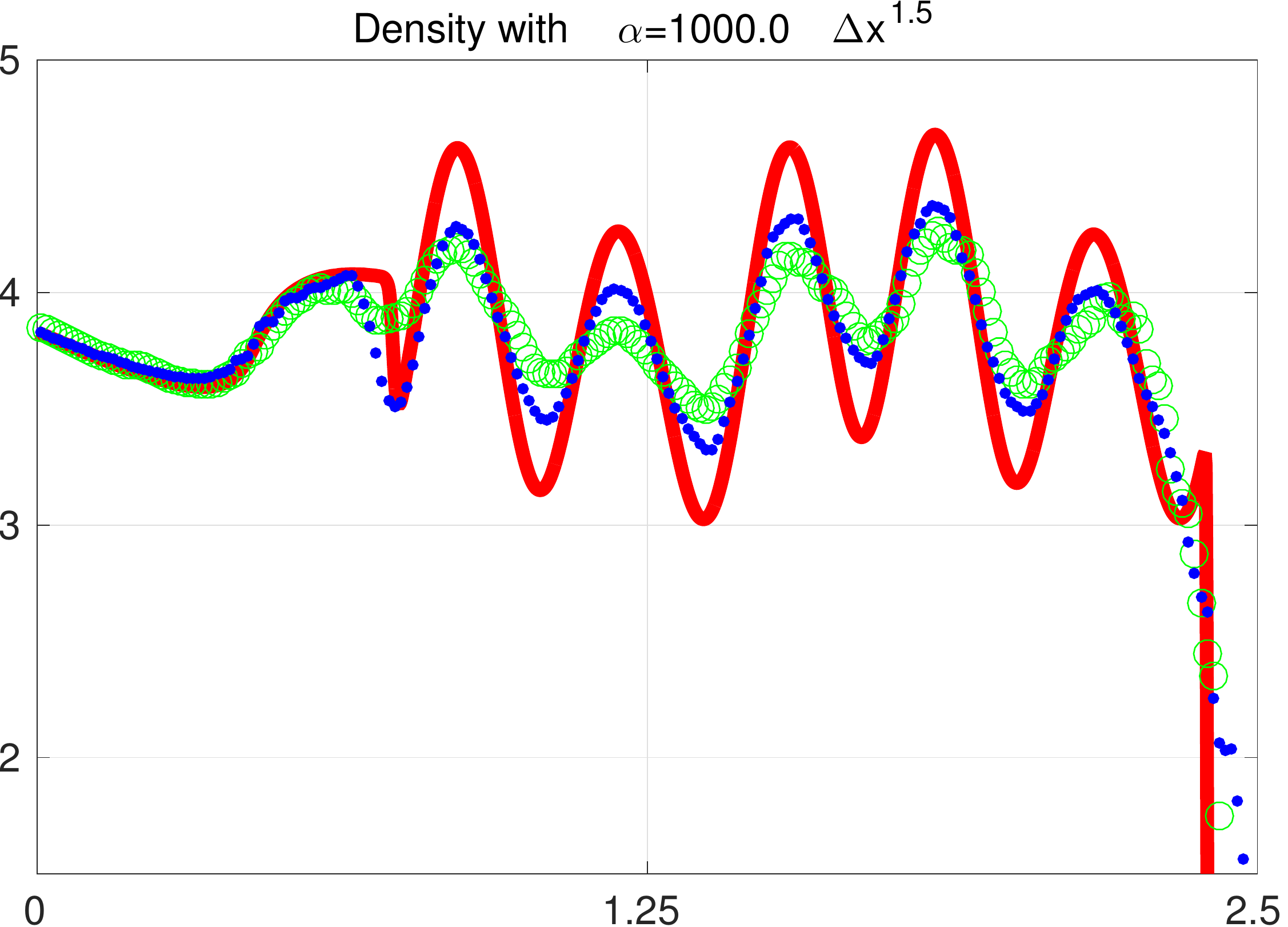}}
 \subfigure[]{\includegraphics[width=0.45 \linewidth]{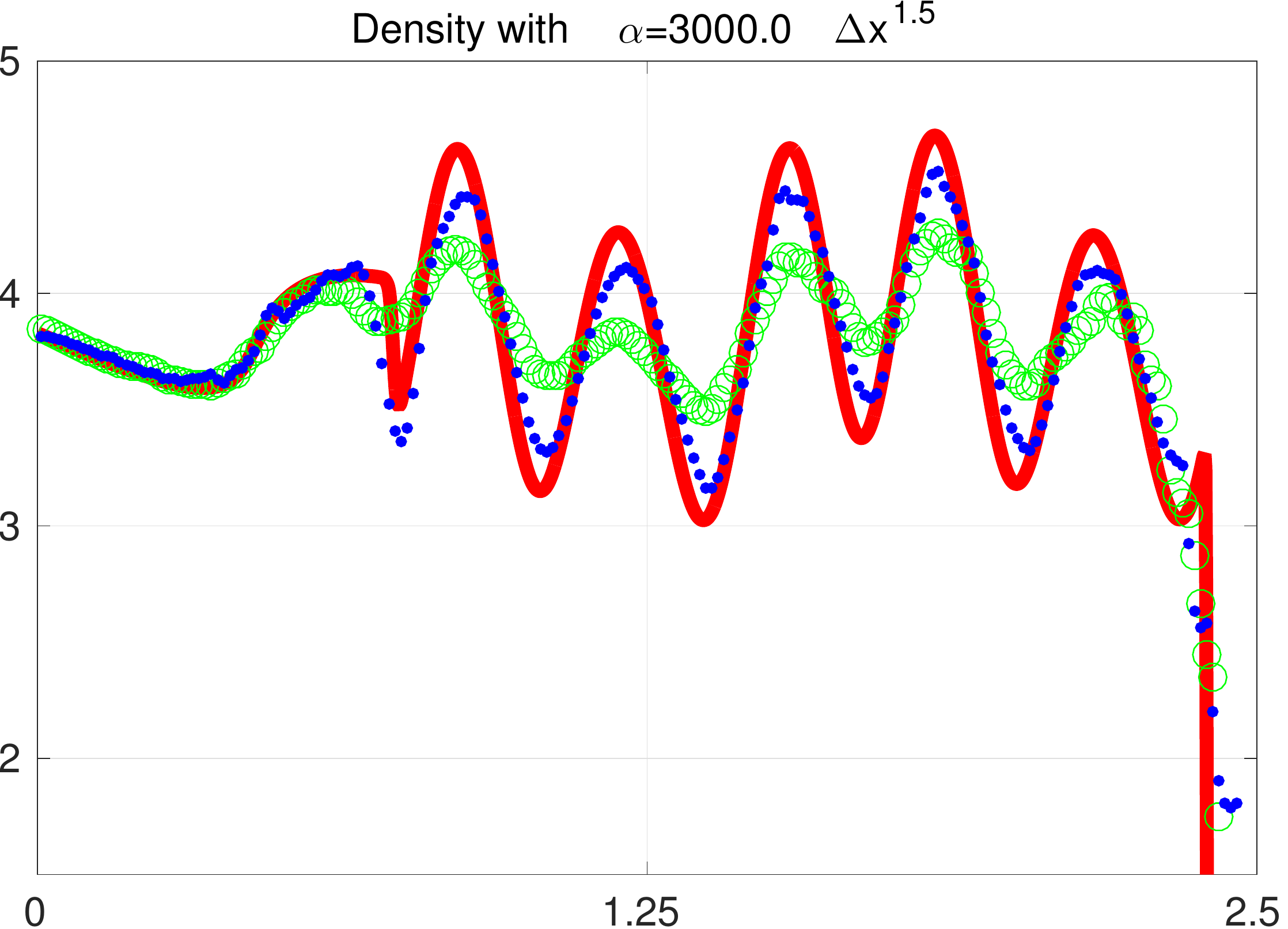}}
  \subfigure[]{\includegraphics[width=0.45 \linewidth]{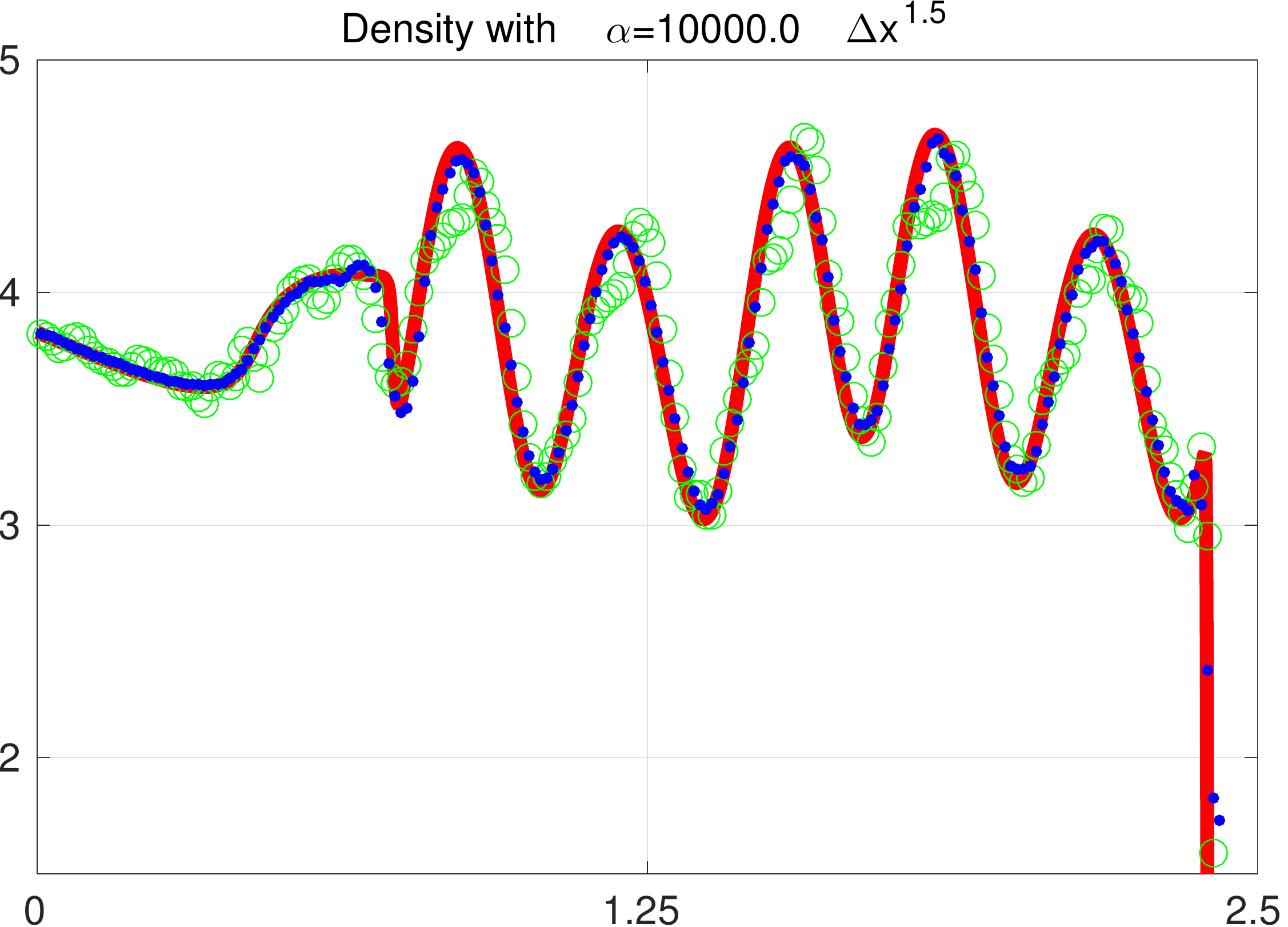}}
  \centering\subfigure{\includegraphics[width=0.45 \linewidth]{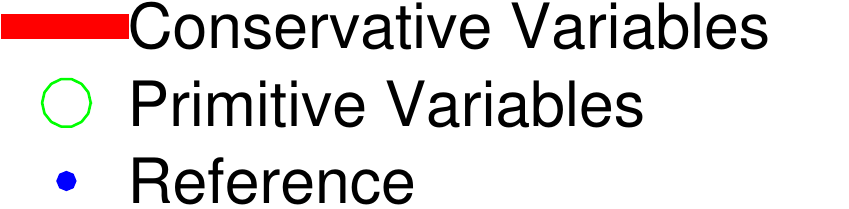}}
 \caption{Shock-entropy problem with $4$ points plotted per cell.  In these
 panels we highlight the benefit of using the primitive 
 variables to conduct the limiting at negligible additional cost.
 To highlight the study, we compare different results that are constructing by
 modifying the subcell tolerance parameter $\alpha$. 
 We observe that selecting a relatively large value of $\alpha$ combined with
 primitive variables is our best choice, given that the conserved variables
 tend to introduce additional oscillations post shock in smooth regimes. 
 \label{fig:pickalpha}
}
 
 \end{figure}

\begin{figure}
\begin{center}
 \includegraphics[width=0.45 \linewidth]{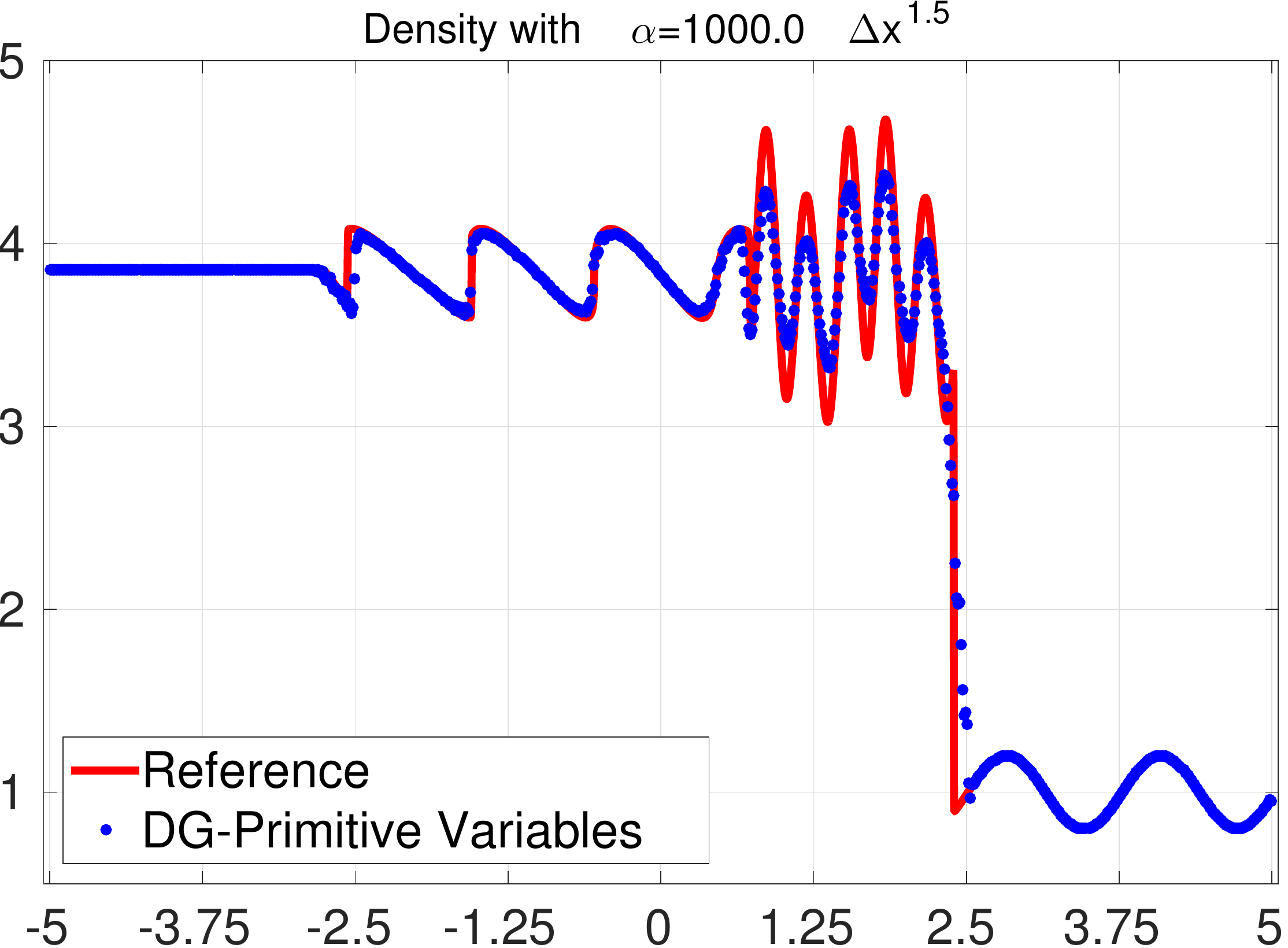}
 \caption{The full plot of the shock-entropy wave problem with $\alpha=500 h^{1.5}$ and primitive variable limiting. This $\alpha$ value
 seems to give a good compromise of resolution and limiting. We will use $\alpha=500 
 h^{1.5}$ in our future calculations unless otherwise noted.}
\label{fig:shuosher}
\end{center}
 \end{figure}

The results of two sets of simulations, one where $\theta_i$ is computed from
sampling conserved variables and one from sampling primitive variables,
is shown in Figure \ref{fig:pickalpha}. For each version of the limiter,
 we highlight the effect of the
choice of $\alpha$ that permits subcell oscillations to persist.  
Additionally, we show the full results of this problem with the primitive
variables used to compute $\theta_i$ and $\alpha = 500 h^{1.5}$
in Figure \ref{fig:shuosher}.
These results clearly show the advantage of using primitive variables in
defining $\theta_i$ and give a range of $\alpha$ values that seem to give
good results.

In every simulation past this point we use primitive variables to define the
rescaling parameter and we take $\alpha = 500 h^{1.5}$, which
seems to balance the need to suppress unphysical oscillations and
still retain high-order accuracy in smooth regions. The only exception to this is the
example in Section \ref{subsec:mach3wind}, where for most of the
computational domain we do use the primitive variables to calculate $\theta_i$,
but there is a small region around the corner of the backward-facing step
where we switch to a more aggressive limiter based on conservative variables. 



\subsection{On the importance of $\alpha$: A simple 2D case study}
\label{subsec:2deuler}
The convergence properties shown in Section \ref{subsec:1dadvection} should still be observable in 2D and for systems. To verify 
this we consider a 2D Euler example that has a smooth exact solution.
This example involves an ideal gas with $\gamma=1.4$, on $(x,y) \in [0,2] \times [0,2]$ with double-periodic boundary conditions
and the following initial conditions:
\begin{equation}
   \left(\rho, u^1, u^2, p\right) = \bigl( 1+0.2 \sin(\pi(x+y)), 0.7, 0.3, 1.0 \bigr).
\end{equation}
In the exact solution $u^1$, $u^2$, and $p$ remain constant for all space and time,
while the density, $\rho$, is advected with the constant fluid velocity given
by $(u_1, u_2) = (0.7, 0.3)$. The numerical $L_2$ errors after one revolution
(i.e., time $t=2$) are shown in Table \ref{table:alpha_convergence2D}. To produce these results we used the same
time-stepping scheme as in Section \ref{subsec:1dadvection}. 
We again observe that we must
have a nonzero value of $\alpha$ to asymptotically match the unlimited solution's convergence even for this very simple test problem.

\begin{table}
 \begin{center}
\hspace*{-50pt}\begin{tabular}{|r||c|c||c|c||c|c|}
\hline
\bf{Mesh} & \bf{{No Limiter}} & \bf{{Order}}& \bf{{$\alpha=0$}} & \bf{{Order}}& \bf{{$\alpha=0.5 \hgrid^{1.5}$}} & \bf{{Order}}\\
\hline
\hline
$   10$ & $1.81\times 10^{-05}$ & --- & $4.24\times 10^{-02}$ & --- & $1.07\times 10^{-02}$ & ---\\
\hline
$   15$ & $3.39\times 10^{-06}$ & $4.97$ & $3.12\times 10^{-02}$ & $0.92$ & $6.07\times 10^{-03}$ & $1.68$ \\
\hline
$  22$ & $6.87\times 10^{-07}$ & $4.48$ & $2.23\times 10^{-02}$ & $0.94$ & $2.12\times 10^{-03}$ & $2.95$ \\
\hline
$  33$ & $1.33\times 10^{-07}$ & $4.05$ & $1.75\times 10^{-02}$ & $0.59$ & $8.56\times 10^{-05}$ & $7.91$ \\
\hline
$  50$ & $2.50\times 10^{-08}$ & $4.36$ & $1.12\times 10^{-02}$ & $1.17$ & $2.50\times 10^{-08}$ & $21.25$ \\
\hline
$  75$ & $4.91\times 10^{-09}$ & $4.02$ & $7.84\times 10^{-03}$ & $0.88$ & $4.91\times 10^{-09}$ & $4.02$ \\
\hline
$  113$ & $9.52\times 10^{-10}$ & $4.15$ & $5.11\times 10^{-03}$ & $1.08$ & $9.52\times 10^{-10}$ & $4.15$ \\
\hline
\end{tabular}\hspace*{-50pt}
  
\caption{ 
    Convergence for the 2D Euler equation problem from Section \ref{subsec:2deuler}. All errors are $L_2$ norm errors. We
    see that with a nonzero value of $\alpha$, we eventually match the unlimited problem's convergence rate, but that without a nonzero
    value of $\alpha$, the convergence is very poor. In fact when the other methods are showing their theoretical asymptotic convergence rates the
    method using $\alpha=0$ is still showing first order convergence (it should be able to achieve second order convergence).}
    
\label{table:alpha_convergence2D}   
 \end{center}
\end{table}

\section{Positivity preservation}
\label{sec:positivity-preservation}
Given the framework of the limiter, the inclusion of the recent and popular
positivity-preserving methods developed by Zhang and Shu \cite{ZhangShu11}
is essentially automatic.  For example, to guarantee positivity of a scalar
function that has a positive mean, we simply replace the computation of the
rescaling parameter in equation \eqref{eqn:thetai} of {\bf Step 4} of the 
algorithm with an additional term to take the minimum over:
\begin{equation*}
     \theta_{i} = \min\left\{ 1, \ \theta_{m_i},\ \theta_{M_i}, \,
     \frac{\bar{q}_i}{ \bar{q}_i - q_{m_i} }  \right\}.
\end{equation*}
Extensions to the compressible 
Euler equations to maintain positive pressure require expanding
the unknown variables as a quadratic function of $\theta$. The details of this can be
found in the work of Zhang and Shu \cite{ZhangShu11}.

\section{Additional numerical results}
\label{sec:numerical-results}


In this section we apply our numerical scheme to several additional standard numerical test cases.
We include one and two-dimensional test cases on both Cartesian and unstructured grids. 
These benchmark test cases are used to verify the accuracy and robustness of the proposed method.



\subsection{2D Riemann problems}

First, we consider two multidimensional
Riemann problem test cases. These problems have been used extensively as
benchmark test cases for other limiters \cite{kurganov2002solution,dumbser2014posteriori}. 

The first example, ${\bf RP1}$, is a well known example that can be found in \cite{MR1241592} and \cite{leveque1997wave}.
We divide the domain $[-0.5,0.5]\times[-0.5,0.5]$ into four rectangular quadrants, all of which meet at the point
$(0.3,0.3)$. We number the quadrants starting in the upper left hand corner, and
continue in a counterclockwise manner.  For example, the upper left hand corner is
Quadrant 1, and the upper right hand corner is Quadrant $4$.
The initial conditions for {\bf RP1} are defined in Table \ref{table:rp1}.

\begin{table}[!t] 

\begin{center}
\begin{tabular}{|c||||c|c|c|c||||c|c|c|c|}
\hline
{{\bf Quadrant}}& $\rho$ & {\bf $u^1$} & {\bf $u^2$} & $p$ & $\rho$ & {\bf $u^1$} & {\bf $u^2$} & $p$ \\
 \hline\hline
 {1} & $0.5323$ & $1.2060$ & $0.0000$ & $0.3000$ & $2.00$ & $0.75$ & $0.50$ & $1.00$ \\
 \hline
 {2} & $0.1380$ & $1.2060$ & $1.2060$&  $0.0290$ & $1.00$ & $-0.75$ & $0.50$&  $1.00$ \\
 \hline
 {3} & $0.5323$ & $0.0000$ & $1.2060$ & $0.3000$ & $3.00$ & $-0.75$ & $-0.50$ & $1.00$ \\
 \hline
 {4} & $1.5000$ & $0.0000$ & $0.0000$ & $1.5000$ & $1.00$ & $0.75$ & $-0.50$ & $0.00$ \\
 \hline
 \end{tabular}
 \end{center}
 \caption{Initial conditions for {\bf RP1} (columns 2--5) and {\bf RP3} (columns 6--9).
 \label{table:rp1}
 }
\end{table}

The second two-dimensional Riemann problem we consider is called {\bf RP3}
\cite{kurganov2002solution,dumbser2014posteriori}.
The domain $[-0.5,0.5]\times[-0.5,0.5]$ is
divided into four equal area quadrants, and the numbering of the quadrants is
also defined in a counterclockwise fashion starting at the upper left hand corner.
The initial conditions for {\bf RP3} are defined in Table \ref{table:rp1}.

\begin{figure}
 \begin{center}
  \subfigure[Density for {\bf RP1} at $t=0.8$. One point per cell.]{\includegraphics[width=0.4\linewidth]{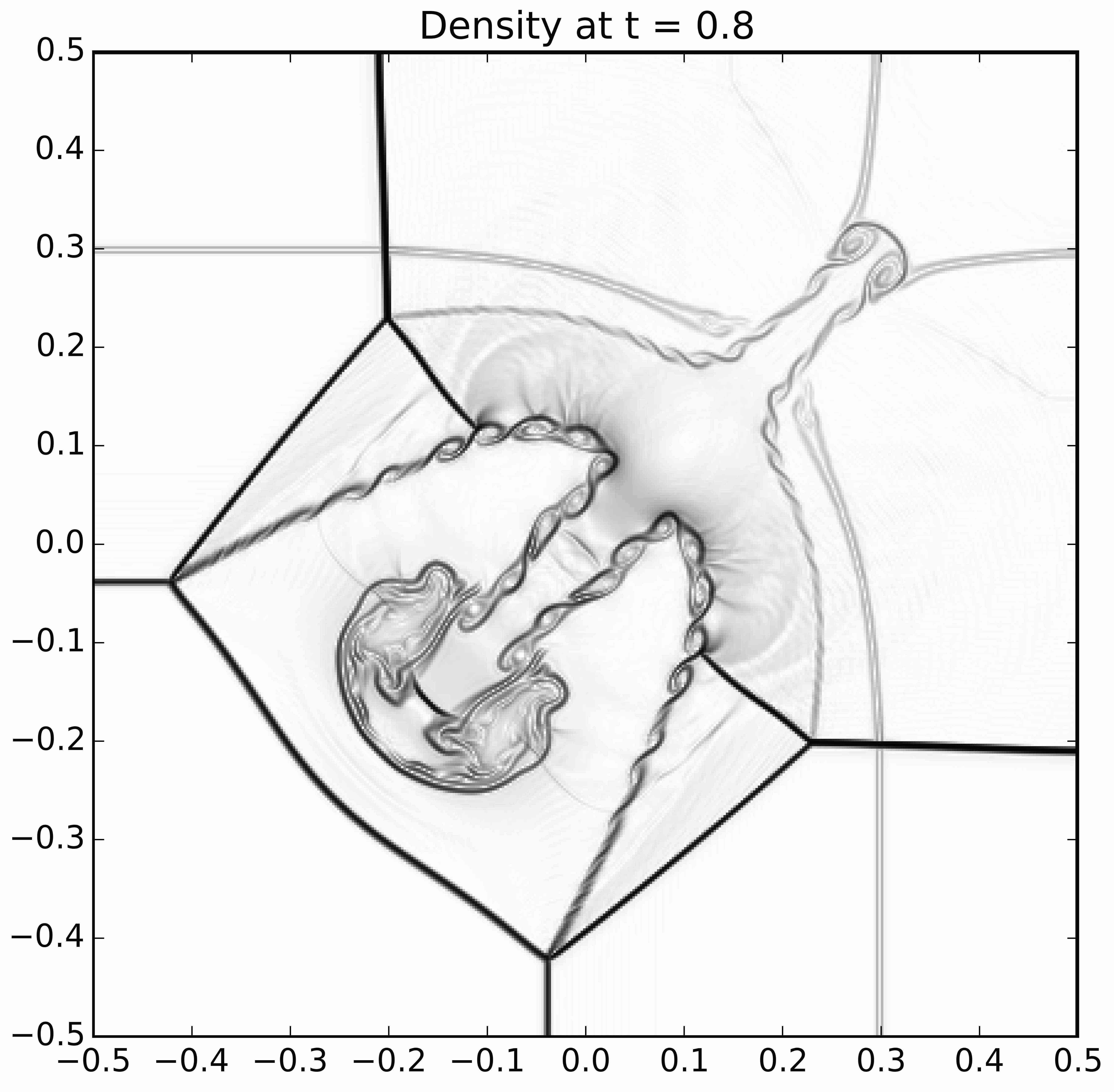}}
  \subfigure[Density for {\bf RP3} at $t=0.3$. One point per cell.]{\includegraphics[width=0.4\linewidth]{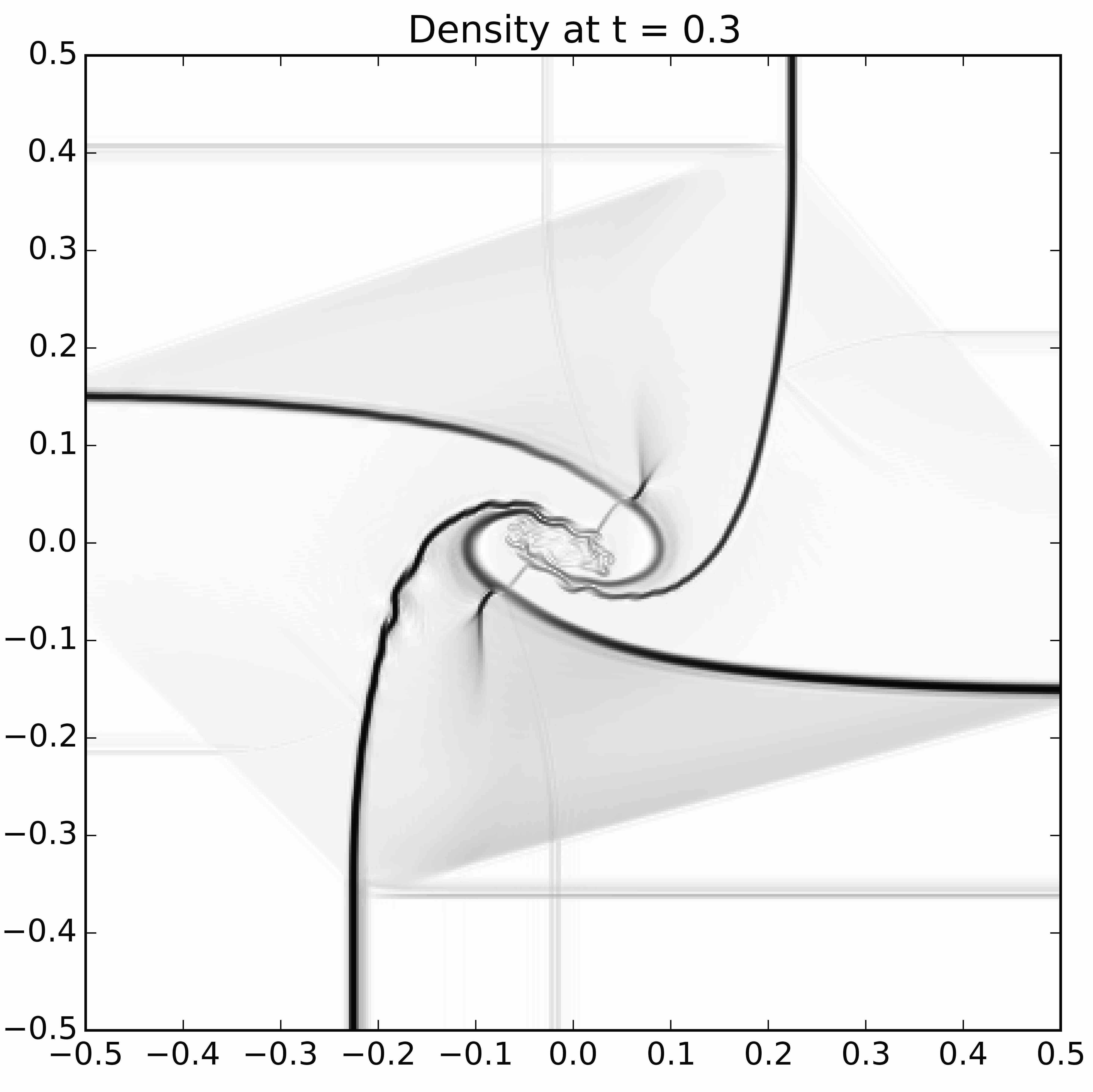}}
  \subfigure[Density for {\bf RP1} at $t=0.8$. Nine points per cell.]{\includegraphics[width=0.4\linewidth]{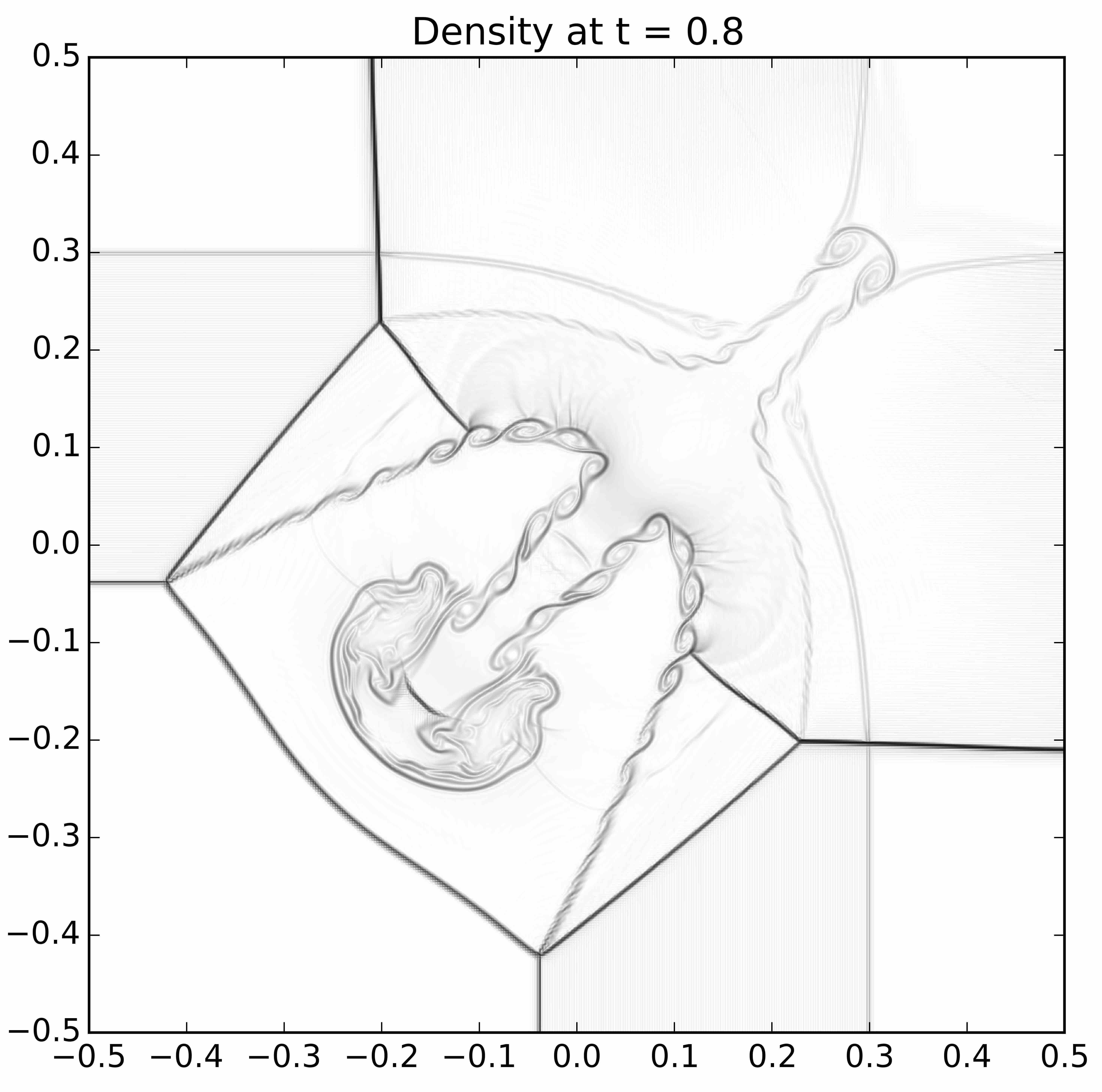}}
  \subfigure[Density for {\bf RP3} at $t=0.3$. Nine points per cell.]{\includegraphics[width=0.4\linewidth]{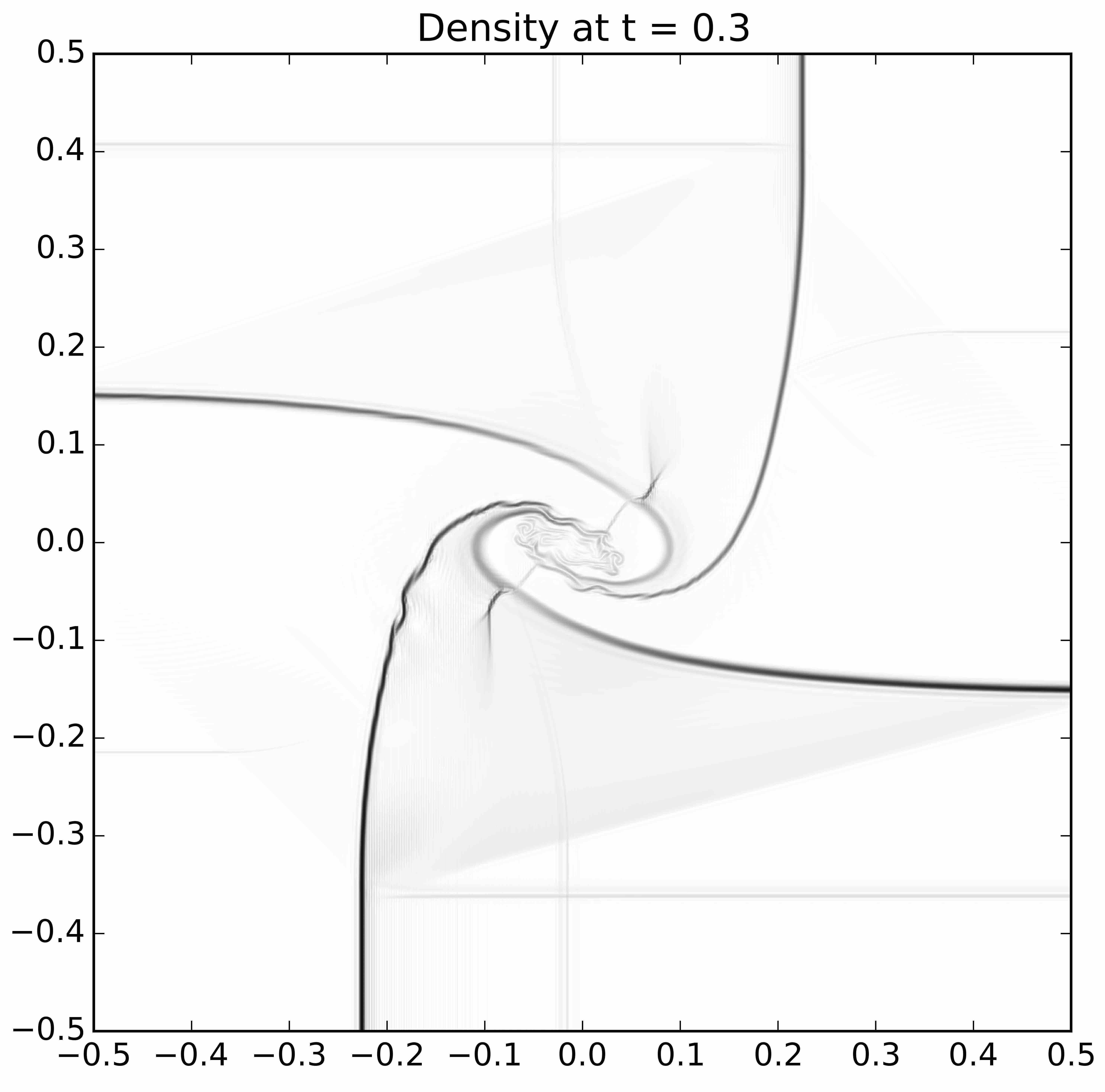}}   
\caption{2D Riemann problems. For {\bf RP1} $30$ contours equally spaced between $0.3$ and $1.8$ are plotted.
For {\bf RP3} $30$ contours equally spaced between $0.1$ and $3.3$ are plotted.}
\label{fig:2dRiemann}
\end{center}
\end{figure}

We run both of these examples on a grid with $400 \times 400$ elements and use
third-order ${\mathbb P}_2$ elements. 
Time-stepping
was done with the classical third-order SSPRK3 method using a CFL number of $0.1$. Both examples are run using
a Rusanov (local Lax Friedrichs (LLF)) flux \cite{article:Ru61}. 
Results are plotted in Figure \ref{fig:2dRiemann}, which are schlieren
plots that show the magnitude of the
gradient of the solution density $|\nabla \rho|$. We show two figures, one with 1 point per cell and the other with three points per direction
per cell. We make this comparison because the limiter permits tiny oscillations to remain
in the solution, since it ignores
all variations smaller than a certain threshold (controlled by the parameter $\alpha$). 
As the mesh is refined, $\alpha$ also shrinks, and therefore the oscillations
will vanish with mesh refinement.  We note that these oscillations are not visible
 in the plots showing one point per cell. 

\subsection{Double Mach reflection}
For our next two-dimensional example, we run the classical double Mach reflection test problem (e.g., see Qiu and Shu \cite{JQiuShu05}). 
In this test problem a Mach $10$ shock approaches a wedge with an angle of $30$ degrees.
In order to simplify the implementation of this problem on a Cartesian grid, it is typical to rotate the whole problem so that the wedge conforms
to the bottom boundary.  
The problem is run on the domain $[0,3.2]\times[0,1]$. For the bottom boundary
we impose the exact solution for $x<\frac{1}{6}$ and we enforce 
solid wall boundary conditions for $x \ge \frac{1}{6}$. 
The initial conditions are
\[
\left( \rho, u^1, u^2, p \right)
=
\begin{cases}
\left( 8.0, 8.25, \cos\left( \frac{\pi}{6} \right), 8.25 \sin\left( \frac{\pi}{6} \right), 116.5 \right)
&  \text{if} \quad x < \frac{1}{6}, \\
\left( 1.4, 0, 0, 1 \right) & \text{otherwise.}
\end{cases}
\]

In the Cartesian case we run this test using a $960 \times 240$ cell grid
matching the second highest resolution in\cite{JQiuShu05}. In the unstructured grid case we use an example with $204,479$ cells. Both
examples are run with the 10-stage SSPRK4 method with a CFL constant of $0.08$.
This smaller CFL number was chosen to guarantee 
positivity of the solution.

\begin{figure}[!ht]
 \begin{center}
  \subfigure[Cartesian grid]{\includegraphics[width=0.95\linewidth]{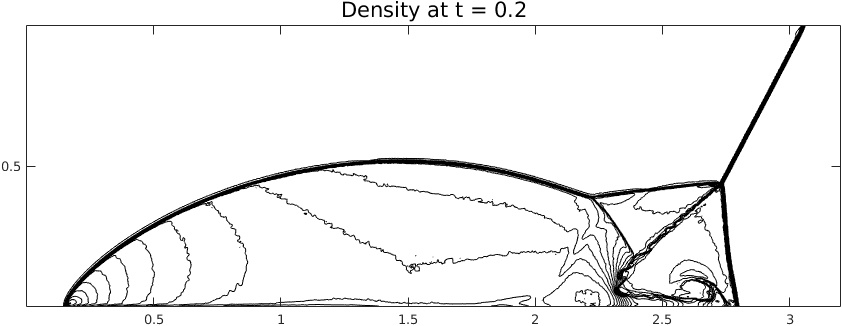}}
  \subfigure[Unstructured grid]{\includegraphics[width=0.95\linewidth]{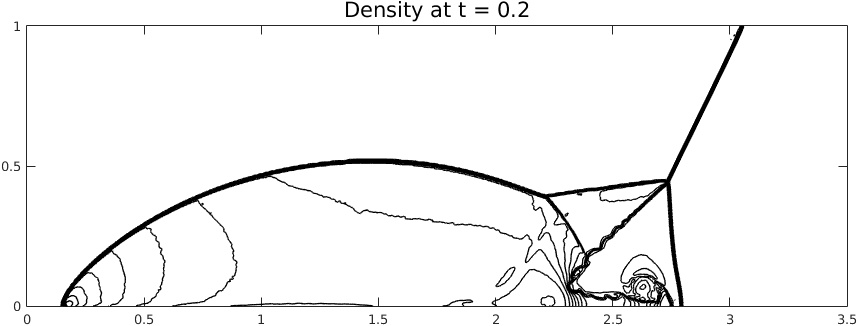}}   
\caption{Double Mach reflection.  
For this example we used an $\alpha$ equal to 
$1000 \hgrid^{1.5}$. We plot a total of $30$ contours that range equispaced from $1.5$ to $22.7$. Both
plots were done using one point per cell.
\label{fig:doublemach}
}
\end{center}
\end{figure}

\begin{figure}[!ht]
 \begin{center}
   \subfigure[$\mathbb{P}_3$ on a Cartesian grid]{\includegraphics[width=0.48\linewidth]{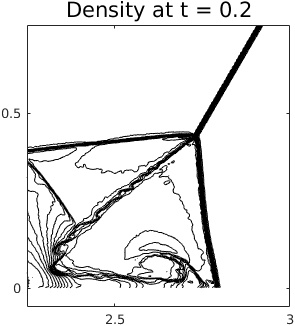}}  
   \subfigure[$\mathbb{P}_3$ unstructured grid]{\includegraphics[width=0.48\linewidth]{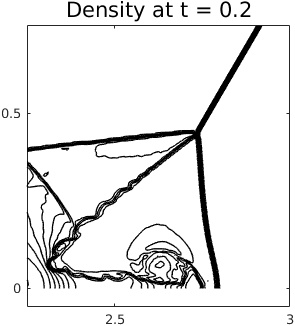}}    
\caption{The output from Figure \ref{fig:doublemach} zoomed in to show the resolution near the contact discontinuity. This contact
discontinuity displays the roll up behavior that is expected in high resolution numerical simulations}

\label{fig:doublemachZoom}
\end{center}
\end{figure}

This problem is challenging given that there is a complicated strong-shock structure, as well as a region where there
is vortical flow. It can be seen in Figure \ref{fig:doublemach} that both limiters essentially succeed in capturing this behavior.
The unstructured grid example appears much better than the Cartesian grid example; this is
 be due to the fact that the unstructured
grid example makes use of all cells touching a current cell, and on the Cartesian grid we only use cells sharing an edge. This was done because 
ideally we do not want to introduce any extra communication in our limiter (by communicating across corners). However in the unstructured
grid case just limiting using cells sharing edges is not good enough and so the unstructured grid limiter had to make use of values in cells
sharing nodes instead of edges. In Figure \ref{fig:doublemachZoom} we see a close of up the solution in the vortical flow region. Both
solutions are certainly forming vortices but the unstructured grid again seems to give much better behavior. However neither example quite matches 
the maximum rollup seen in the comparable resolution examples in \cite{JQiuShu05}. 

\subsection{Mach 3 wind tunnel with a step on Cartesian and unstructured meshes}
\label{subsec:mach3wind}
We next consider a classical test problem \cite{woodward1984numerical,JiangShu96} that
is often overlooked, likely due to its difficulty.
The problem involves an ideal gas with $\gamma=1.4$ and a strong shock moving at Mach $3$
that passes past a forward facing step. The initial conditions for this
problem are given by
$
\left( \rho, u^1, u^2, p \right)
  =
\left( 1.0, 3.0, 0.0, {\gamma}^{-1} \right).
$
The computational domain is the rectangle $[0,3]\times[0,1]$
with the region $[0.6,3.0]\times[0.0,0.2]$ cut out. In the Cartesian grid
case, we use a uniform grid with $480 \times 160$ elements,
and we use a total of $90,342$ cells for the unstructured grid case.
The Cartesian grid matches
the moderately refined implementation in \cite{CoShu98}, however we run with
the higher-order $\mathbb{P}_3$ elements instead of $\mathbb{P}_2$ elements
used in \cite{CoShu98}. 
Both examples are run with SSPRK4 using the LLF flux and a CFL of $0.4$.

This is a complicated problem featuring a long contact discontinuity where,
again, vortices that should form are observable with high-resolution simulations. 
An additional and well-known difficulty is that this example tends to form a spurious entropy
layer above the step. This results in the creation of a stem-like structure where
there would normally be a shock reflecting off of the step. This is especially noticeable
in the Cartesian grid plot. To mitigate
this effect, we use a slightly more aggressive implementation of our
limiter near top of the step.  In place of using point-wise values of
the primitive variables to
estimate local upper and lower bounds, we use cell averages of the conserved
variables in this edge region only.
That is, we replace $M^{\ell}_i$ and $m^{\ell}_i$ in 
equations \eqref{eqn:big_mi_sys} and \eqref{eqn:little_mi_sys} by 
\begin{equation*}
     M_i := \max \left\{ \bar{q}_i+\alpha(\hgrid), \, \max_{j \in N_{\Tm_i}}\left\{\bar{q}^h_{M_j}
     \right\}
     \right\} \, \text{ and } \,
     m_i := \min\left\{\bar{q}_i-\alpha(\hgrid), \, \min_{j \in N_{\Tm_i}}\left\{ \bar{q}^h_{m_j}
     \right\}
     \right\},
\end{equation*}
where the $\ell$ superscript has  been suppressed. 
Note that this means that in this special region we are not using the primitive
variables for determining $\theta_i$, but actually the conserved variables; furthermore,
we are using the average values $\bar{q}$ rather than point values at the
points defined by $\chi_i$.
Additionally we use a much smaller $\alpha$ value in this region to obtain
better results. In the majority of the domain we take 
$\alpha=500 \hgrid^{1.5}$,
but in the small region above the step we take
$\alpha=20 \hgrid ^{1.5}$.
The region using the more aggressive limiter extends $6$ cells above the step in the Cartesian grid case, and until $y=0.24$
above the step in the unstructured
grid case. These values can certainly be modified, but we found that if the region is taken to be only one cell width, then the entropy
layer tends to be much more prominent.  Moreover, if the more aggressively limited region is too large the shocks in the problem
will be smeared out too heavily.

\begin{figure}
 \begin{center}
  \subfigure[Cartesian grid]{\includegraphics[width=0.87\linewidth]{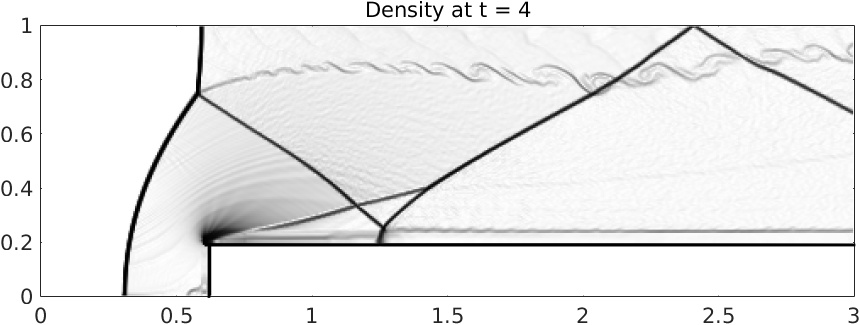}}
   \subfigure[Unstructured grid]{\includegraphics[width=0.87\linewidth]{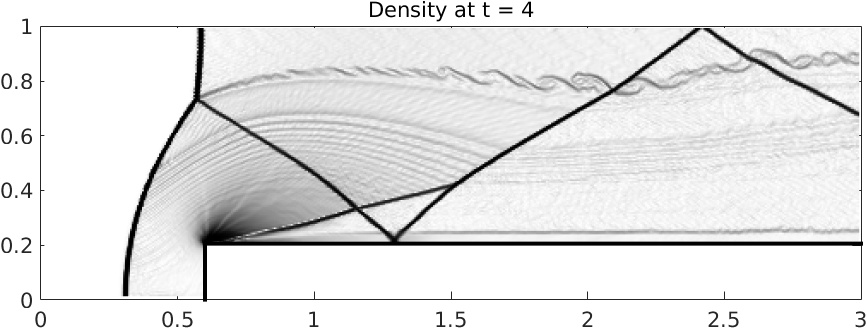}}
\caption{Density schlieren plots from the forward facing step problem run on Cartesian and unstructured grids. The Cartesian grid used
$\hgrid=\frac{1}{160}$, and the unstructured grid used $90,342$ cells (that were uniformly sized). Both were plotted with one point per cell.}
%
\label{fig:WindTunnel}
\end{center}
\end{figure}

\begin{figure}
 \begin{center}
  \subfigure[Cartesian grid]{\includegraphics[width=0.87\linewidth]{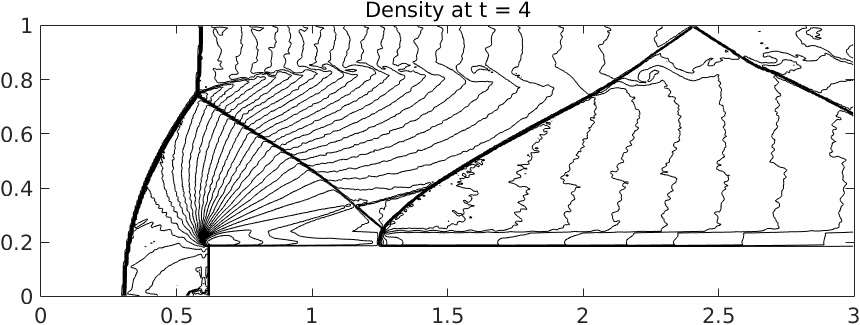}}
   \subfigure[Unstructured grid]{\includegraphics[width=0.87\linewidth]{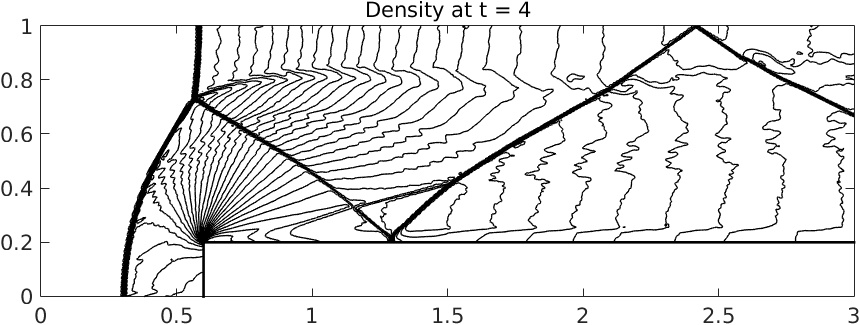}}
\caption{Density contour plots from the forward facing step problem run on Cartesian and unstructured grids. There are $30$ contours
spaced evenly from $0.1$ to $4.54$.}
\label{fig:WindTunnelContour}
\end{center}
\end{figure}

In Figures \ref{fig:WindTunnel} and \ref{fig:WindTunnelContour} we observe that both methods do an excellent job of capturing the structure of the contact 
discontinuity. However it appears that the Cartesian grid case suffers from a pronounced spurious entropy layer; although it should be noted that the entropy
layer is not atypically large (e.g., see \cite{CoShu98,ZhuZhoShuQiu13}). Given the relative
simplicity of this limiter the results seem very promising. 

%
%
%

\section{Conclusions}
\label{sec:conclusions}
In this work, we have presented a novel limiter for the discontinuous Galerkin
method.  It is the first of its kind given its ease of implementation, ability
to retain high-order accuracy, and robustness for simulation problems with strong
shocks. The results are achieved without the need for an additional shock-detector. 
The proposed method can 
be viewed as either an extension of the Barth-Jespersen \cite{BarthJesp89}
limiter to discontinuous Galerkin methods or as extending the Zhang and Shu
\cite{ZhangShu11}
positivity preserving method to be able to accommodate shocks and not just 
satisfy the positivity principle.  
Several numerical results were presented, including single and
multidimensional examples on Cartesian and unstructured grids.  Future work
will involve investigating the ability of this limiter to work on more
complicated hyperbolic problems such as magnetohydrodynamics (MHD) and three
dimensional problems with structured and unstructured grids.

\bigskip

\noindent
{\bf Acknowledgements.}
The work of SAM was supported in part by NSF grant DMS--1216732. 
The work of JAR was supported in part by NSF grant DMS--1419020.

\bibliographystyle{plain}
\bibliography{References}

\end{document}